\documentclass[final]{siamltex}
\usepackage{amsfonts,amsmath,mathrsfs,bm}
\usepackage{amssymb}
\usepackage[notref,notcite]{showkeys}
\usepackage{hyperref}
\usepackage{graphicx}
\usepackage{wrapfig}
\usepackage{subfig}
\usepackage{float}
\usepackage{bbm}
\usepackage{algorithmic, algorithm}
\usepackage{gensymb,color,soul}
\usepackage{tikz}
\usepackage{tikz-3dplot}

\newcommand{\N}{\mathbb{N}}
\newcommand{\R}{\mathbb{R}}
\newcommand{\C}{\mathbb{C}}

\newtheorem{remark}[theorem]{Remark}

\DeclareMathOperator{\argmin}{\rm arg\, min}
\DeclareMathOperator{\argmax}{\rm arg\, max}

\title{Bayesian experimental design for head imaging by electrical impedance tomography}

\author{
N.~Hyv\"onen\footnotemark[2]
\and A.~J\"a\"askel\"ainen\footnotemark[2]
\and R.~Maity\footnotemark[2]
\and A.~Vavilov\footnotemark[2]
}

\begin{document}
\maketitle

\renewcommand{\thefootnote}{\fnsymbol{footnote}}
\footnotetext[2]{Aalto University, Department of Mathematics and Systems Analysis, P.O.~Box 11100, FI-00076 Aalto, Finland (nuutti.hyvonen@aalto.fi, altti.jaaskelainen@aalto.fi, ruma.maity@aalto.fi, anton.vavilov@aalto.fi). This work was supported the Jane and Aatos Erkko Foundation and the Academy of Finland (decision 348503, 353081).}

\begin{abstract}
This work considers the optimization of electrode positions in head imaging by electrical impedance tomography. The study is motivated by maximizing the sensitivity of electrode measurements to conductivity changes when monitoring the condition of a stroke patient, which justifies adopting a linearized version of the complete electrode model as the forward model. The algorithm is based on finding a (locally) A-optimal measurement configuration via gradient descent with respect to the electrode positions. The efficient computation of the needed derivatives of the complete electrode model is one of the focal points. Two algorithms are introduced and numerically tested on a three-layer head model. The first one assumes a region of interest and a Gaussian prior for the conductivity in the brain, and it can be run offline, i.e., prior to taking any measurements. The second algorithm first computes a reconstruction of the conductivity anomaly caused by the stroke with an initial electrode configuration by combining lagged diffusivity iteration with sequential linearizations, which can be interpreted to produce an approximate Gaussian probability density for the conductivity perturbation. It then resorts to the first algorithm to find new, more informative positions for the available electrodes with the constructed density as the prior.
\end{abstract}

\renewcommand{\thefootnote}{\arabic{footnote}}



\pagestyle{myheadings}
\thispagestyle{plain}
\markboth{N.~HYV\"ONEN, A.~J\"A\"ASKEL\"AINEN, R.~MAITY, AND A.~VAVILOV}{EXPERIMENTAL DESIGN FOR HEAD IMAGING BY EIT}

\section{Introduction}
\label{sec:introduction}

This work considers Bayesian {\em optimal experimental design} (OED) for choosing positions of electrodes for monitoring stroke by {\em electrical impedance tomography} (EIT). The goal of EIT is to reconstruct (useful information on) the internal conductivity of an imaged physical body from boundary measurements of current and voltage at contact electrodes attached to the boundary of the body. This task constitutes an inverse boundary value problem that is nonlinear and severely illposed. For general information on EIT, we refer to the review articles \cite{Borcea02,Cheney99,Uhlmann09}.

We resort to a linearized version of the smoothened {\em complete electrode model} (CEM) \cite{Cheng89,Somersalo92, Hyvonen17b} as the forward model for the EIT measurements and aim to choose the positions of the available finite number of electrodes so that they maximize the information on the conductivity changes in a {\em region of interest} (ROI) that corresponds to an area of the brain affected by a stroke. Our study is motivated by the potential application of EIT to bed-side monitoring of stroke in intensive care~\cite{Toivanen21}: When a stroke patient enters a hospital, a {\em computed tomography} (CT) image of the patient's head is taken, which reveals the head anatomy and the location of the stroke. In particular, the linearization of the forward model around a background conductivity can be motivated by the goal of detecting regrowth of a hemorrhagic stroke or secondary bleeding induced by an ischemic stroke, both of which are presumed to cause moderate conductivity changes around the initially damaged part of the brain. However, the linearization also considerably simplifies the application of Bayesian OED as explained in the following.

A Bayesian optimal design $d^*$ is defined as an admissible maximizer for the expectation of a given utility function $\Upsilon(w, y; d)$, with the expectation taken with respect to the data $y \in \mathcal{Y}$ and the parameter of interest $w \in \mathcal{W}$ \cite{alexanderian2021optimal_review, chaloner1995bayesian, rainforth2023modern, ryan2016review}. In mathematical terms,
\begin{equation}
	\label{eq:OED_task}
	d^* 
	 =  \underset{d\in \mathcal{D}}{\argmax} \int_{\mathcal{Y}} \int_\mathcal{W} \Upsilon(w, y; d) \pi(w, y \, | \,  d) \, {\rm d}w \, {\rm d}y,
\end{equation}
where $\mathcal{D}$ is the space of admissible designs and $\pi(w, y \, | \,  d)$ is the joint probability density of the parameter and data for the design $d$. In this work, $\Upsilon$ is chosen to be a \emph{negative quadratic loss function} that measures distance from $w$ to the mean $\widehat{w}(y; d)$ of the posterior density $\pi(w \, | \, y; d)$, which leads to the concept of Bayesian A-optimality in our linear(ized) setting. Another common choice for $\Upsilon$ is related to the Kullback--Leibler distance between the posterior and prior distributions, which would in turn correspond to Bayesian D-optimality. D-optimality is not considered in this work, although it could be tackled with the same techniques,~cf.~\cite{Hyvonen14}.

For the considered inverse problem of EIT, $d$ determines the positions of the electrodes involved in the design process, $w$ represents the (discretized) change in the conductivity of the imaged head, and $y$ carries the measured potential differences at the electrodes. Due to the employed linearization of the relation between $w$ and $y$, a Gaussian prior for $w$ and an additive Gaussian noise model that is independent of $w$ enable explicit evaluation of the double integral in \eqref{eq:OED_task}, leading to minimization of a weighted trace of the posterior covariance with respect to $d$. As mentioned above, this corresponds to finding an A-optimal design~\cite{alexanderian2016bayesian,chaloner1995bayesian}. The posterior covariance allows a simple representation in terms of the linearized forward map and the prior and noise covariances, which enables performing the experimental design offline assuming the additive noise process is independent of $d$. Take note, however, that minimizing the trace of the posterior covariance with respect to $d$ by a gradient-based minimization algorithm requires evaluating derivatives of the linearized forward model with respect to the electrode positions,~i.e.,~computing certain {\em second} order derivatives of the forward solutions. This setting of experimental design for EIT has been considered in two spatial dimensions in \cite{Hyvonen14} for the standard CEM, but the computation of the required derivatives is computationally more demanding in our three-dimensional setting than in~\cite{Hyvonen14}.

In addition to the offline optimization of the positions of the electrodes, we also consider an alternative approach that is not as straightforwardly applicable to the bed-side monitoring task used as a motivation above due to the requirement to be able to alter the positions of the electrodes attached to the head of the patient. However, it has potential to not require a CT image of the patient, cf.~\cite{Candiani21}. The idea is to assume an edge-promoting (weighted and smoothened) {\em total variation} (TV) prior~\cite{Rudin92} for the conductivity change in the ROI and to choose an initial set of electrodes for which measurement data is measured (or, in our studies, simulated). A preliminary reconstruction is then formed by combining sequential linearizations and the lagged diffusivity iteration~\cite{Vogel96}; see~\cite{Harhanen15} for the introduction of this approach to EIT and~\cite{Candiani21} for its application to head imaging without precise knowledge of the patient's anatomy. In particular, the initial reconstruction accounts for the nonlinearity of the forward model. According to a Bayesian interpretation \cite{Bardsley18,Calvetti08,Helin22}, the lagged diffusivity iteration also produces a Gaussian approximation for the non-Gaussian posterior induced by the TV prior. The covariance of this approximation can then be employed to find more informative positions for the available electrodes in the same manner as in the offline approach described above, with the aim of using the new positions to better monitor possible further changes in the conductivity close to the initially detected one,~i.e.,~close to the region initially affected by the stroke. A closely related technique for dealing with a TV prior in Bayesian OED has previously been considered for the linear inverse problems of X-ray and magnetorelaxometry imaging in \cite{Helin23,Helin22}.

Our approach to Bayesian OED is based on linearization and discretization of an inverse boundary value problem for an elliptic partial differential equation. Although convenient in our computationally demanding setting, especially with the TV prior that is not discretization invariant~\cite{Lassas04}, these steps of linearization and discretization are not necessary for formulating and solving Bayesian OED problems for models involving partial differential equations. We refer the reader to~\cite{alexanderian2021optimal_review} for information on the state-of-the-art for such OED problems and note that certain stability properties of Bayesian OED with respect to approximations in the forward model have recently been studied in~\cite{duong2022stability}, providing a first step toward formally justifying our approach. In particular, the optimization of electrode positions in EIT has previously been considered in \cite{Haber10,Hyvonen14,Karimi20,Smyl20}. The papers \cite{Hyvonen14,Karimi20,Smyl20} tackle two-dimensional settings: \cite{Hyvonen14} essentially corresponds to our offline optimization approach in two-dimensions, \cite{Smyl20} applies deep learning to the OED problem, and \cite{Karimi20} works with the less accurate shunt model~\cite{Cheng89} and evaluates the expected utility by resorting to a double-loop Monte Carlo method. In \cite{Haber10}, a large set of {\em point electrodes} (cf.~\cite{Hanke11b}) is set at predefined positions, and an optimization method with sparsity constraints for the current injections and potential measurements is applied to select a feasible set of electrodes.  

This article is organized as follows. The CEM, its linearization and discretization, as well as the needed derivatives with respect to the electrode positions are introduced in Section~\ref{sec:CEM}. Section~\ref{sec:Bayesian} recalls the basics of finite-dimensional Bayesian inversion, including the Bayesian interpretation of the lagged diffusivity iteration and the reconstruction algorithm based on sequential linearizations used in connection to the TV prior. Section~\ref{sec:oed} considers OED and introduces the related numerical optimization algorithms. Section \ref{sec:numerics} documents the numerical experiments, and Section~\ref{sec:conclusion} contains the concluding remarks.

\section{Complete electrode model}
\label{sec:CEM}

This section recalls the CEM and its derivatives with respect to the conductivity and the electrode contacts and positions. We also introduce the discretized computational models that are used in the numerical experiments.

\subsection{Forward problem}
\label{sec:forward_problem}
Let $M \in \N \setminus \{ 1 \}$ contact electrodes $E_1, \dots,  E_M$ be attached to the boundary of the imaged head that is modeled as a bounded Lipschitz domain $\Omega \subset \R^3$. The electrodes are represented as open connected subsets $E_1, \dots, E_M$ of the boundary $\partial \Omega$ with mutually disjoint closures. The union of the electrodes is denoted by $E = \cup E_m$. A single measurement of EIT consists in driving a net current pattern $I \in \C^M_\diamond$ through the electrodes and measuring the resulting constant electrode potentials $U \in \C^M_\diamond$, with
\[ 
\C^M_\diamond = \Big\{V\in\C^M\,\Big|\, \sum_{m=1}^M V_m = 0\Big\}.
\]
The current pattern $I$ belongs to $\C^M_\diamond$ due to conservation of charge and $U$ based on an appropriate choice for the ground level of potential. The latter vector is identified in the following with the piecewise constant function 
\begin{equation}
\label{eq:piecewise}
U \, = \, \sum_{m=1}^M U_m \chi_m,
\end{equation}
where $\chi_m$ is the characteristic function of $E_m$ on $\partial \Omega$. 

A mathematical model that accurately predicts real-life EIT measurements is the CEM~\cite{Cheng89}. Here we employ its smoothened version \cite{Hyvonen17b} that carries better regularity properties without severely compromising modeling accuracy. According to the smoothened CEM, the electromagnetic potential $u$ inside $\Omega$ and the potentials on the electrodes $U$ satisfy the elliptic boundary value
\begin{equation}
\label{eq:cemeqs}
\begin{array}{ll}
\displaystyle{\nabla \cdot(\sigma\nabla u) = 0 \qquad}  &{\rm in}\;\; \Omega, \\[6pt]
{\displaystyle {\nu\cdot\sigma\nabla u} = \zeta (U - u) } \qquad &{\rm on}\;\; \partial \Omega, \\[2pt]
{\displaystyle \int_{E_m}\nu\cdot\sigma\nabla u\,{\rm d}S} = I_m, \qquad & m=1,\ldots,M, \\[4pt]
\end{array}
\end{equation}
interpreted in the weak sense and with $\nu$ denoting the exterior unit normal on $\partial\Omega$. The conductivity $\sigma$ is assumed to be isotropic and belong to
\begin{equation}
\label{eq:sigma}
L^\infty_+(\Omega) := \{ \varsigma \in L^\infty(\Omega) \ | \ {\rm ess} \inf {\rm Re} (\varsigma) > 0 \},
\end{equation}
whereas the contact conductivity $\zeta$ is modeled as an element of
\begin{equation}
  \label{eq:zeta}
\mathcal{Z} := \big\{ \xi \in L^\infty(E) \  \big| \   {\rm Re}\, \xi \geq 0 \   {\rm and} \  {\rm ess} \sup \big( {\rm Re} ( \xi |_{E_m} ) \big) > 0 \  {\rm for} \ {\rm all} \ m= 1, \dots, M \big\},
\end{equation}
which is interpreted as a subset of $L^\infty(\partial \Omega)$ via zero-continuation. The variational formulation of \eqref{eq:cemeqs} is to find $(u,U) \in \mathcal{H}^1 := H^1(\Omega) \oplus \C^M_\diamond$ such that
\begin{equation}
\label{eq:weak}
B_{\sigma,\zeta}\big((u,U),(v,V)\big)  \,=  \, I\cdot V \qquad {\rm for} \ {\rm all} \ (v,V) \in  \mathcal{H}^1,
\end{equation}
where $\cdot$ denotes the real dot product and the bounded and coercive bilinear form $B_{\sigma,\zeta}: \mathcal{H}^1 \times \mathcal{H}^1 \to \C$ is defined by
\begin{equation}
\label{eq:sesqui}
B_{\sigma,\zeta}\big((w,W),(v,V)\big) = \int_\Omega \sigma\nabla w\cdot \nabla v \,{\rm d}x + \int_{\partial \Omega} \zeta (W-w)(V-v)\,{\rm d}S.
\end{equation}
Note that we intentionally exclude complex conjugations and work with bilinear forms instead of sesquilinear ones because such an approach simplifies certain expressions in what follows. Under the above assumptions, the potential pair $(u,U) \in \mathcal{H}^1$ is uniquely determined by \eqref{eq:weak}~\cite{Hyvonen17b}. For physical justification of \eqref{eq:cemeqs}, see \cite{Cheng89,Hyvonen17b}.

The functional dependence of the electrode potentials on the conductivities and the applied electrode currents is denoted as $U(\sigma,\zeta; I) \in \C^M_\diamond$. In what follows, we employ a basis $I_1, \dots, I_{M-1}$ for $\C^M_\diamond$ as the electrode current patterns and adopt the notation
\begin{equation}
\label{eq:forward_map}
\mathcal{U}\big(\sigma,\zeta; I_1, \dots, I_{M-1}\big) 
= \left[U(\sigma,\zeta; I_1)^{\top}, \dots, U(\sigma,\zeta; I_{M-1})^{\top}\right]^{\top}
\in \C^{M(M-1)}.
\end{equation}
As the current basis is typically either clear from the context or unimportant for the presented arguments, the dependence of $\mathcal{U}$ on $I_1, \dots, I_{M-1}$ is often suppressed.

\subsection{Derivatives of the CEM solution}
We next consider the derivatives of the forward solution $(u,U) \in \mathcal{H}^1$ with respect to the conductivities and the electrode positions. For simplicity, let $\partial \Omega$ and $\partial E$ be of class $\mathcal{C}^\infty$ and assume the conductivities in \eqref{eq:cemeqs} satisfy $\sigma \in L^\infty_+(\Omega) \cap \mathcal{C}^{0,1}(\overline{\Omega})$ and $\zeta \in \mathcal{Z} \cap H^1(\partial \Omega)$, although this level of regularity is not needed for all results presented in this section.\footnote{In particular, it would be sufficient for all presented results to only assume that $\sigma$ is Lipschitz continuous in some neighborhood of the boundary $\partial \Omega$.} We refer to \cite{Hyvonen14} for careful analysis on the existence of the considered derivatives in the less regular case that $\zeta$ is piecewise constant and to \cite{Candiani19,Hyvonen17b} for considerations on differentiation with respect to electrode positions when $\zeta \in H^1(\partial \Omega)$. Apart from assuming a smoother contact conductivity, the main difference to \cite{Hyvonen14} is that we deduce more efficient computational formula for the mixed derivative with respect to the internal conductivity and the electrode positions. Since we consider the presented results straightforward, yet arguably tedious modifications of material in \cite{Hyvonen17b,Hyvonen14}, we do not present detailed proofs. Moreover, we do not claim that our regularity assumptions are optimal; the conditions on $\sigma$ and its perturbations are required in \cite{Hyvonen14} due to the discontinuities in $\zeta$ that make the derivatives with respect to the electrode positions singular in a certain sense \cite{Darde12}.

To begin with, let us define how electrodes are moved along $\partial \Omega$. Let $a \in \mathcal{C}^1(E, \R^3)$ and define a perturbed set of electrodes via
\begin{equation}
\label{eq:perturbed}
E_m^a = \big\{ P_x \big(x +a(x) \big) \, \big| \, x \in E_m \big\} \subset \partial \Omega, \qquad m=1,\dots, M,
\end{equation}
where $P_x: \R^3 \supset B_\rho(x) \to \partial \Omega$ is the projection in the direction of $\nu(x)$ onto $\partial \Omega$, and $B_\rho(x)$ is an open ball of small enough radius $\rho > 0$ centered at $x$. The contact conductivity $\zeta: \partial \Omega \to \C$ is assumed to stretch accordingly, that is, the conductivity on the perturbed electrodes is defined via
\begin{equation}
\label{eq:perturbedz}
\zeta^a\big(P_x (x + a(x) ) \big) = \zeta(x), \qquad x \in E,
\end{equation}
and it vanishes on $\partial \Omega \setminus \overline{E^a}$. It can be shown that there exists a constant $c = c(E,\Omega) >0$ such that the definitions \eqref{eq:perturbed} and \eqref{eq:perturbedz} are unambiguous if $\|a \|_{\mathcal{C}^1(E, \R^3)} < c$; see~\cite{Darde12} for the proof of a closely related result. For a fixed set of electrodes $E_1, \dots, E_M$, we can now interpret the electrode potentials (as well as the interior potential) as a function of four variables, that is,
\begin{equation}
  \label{eq:tuned_forward_map}
U:
\left\{
\begin{array}{l}
(\sigma, \zeta, a, I) \mapsto U(\sigma, \zeta, a; I), \\[1mm]
  L^\infty_+(\Omega) \times \mathcal{Z} \times \mathcal{B}_c \times \C^M_\diamond \to \C^M_\diamond,
\end{array}
\right.
\end{equation}
where $U(\sigma, \zeta, a; I)$ is the second part of the solution to \eqref{eq:weak} with $E_1, \dots, E_M$ replaced by $E_1^a, \dots, E_M^a$, and $\mathcal{B}_c$ denotes the origin-centered open ball of radius $c$ in the topology of $\mathcal{C}^1(E, \R^3)$.

Let us then introduce the required derivatives. It is well known that the map $U: L^\infty_+(\Omega) \times \mathcal{Z} \times \mathcal{B}_c \times \C^M_\diamond \to \C^M_\diamond$ is Fr\'echet differentiable with respect to its first variable, but the same also holds for the second variable~\cite{Darde21} although it is a less obvious result due to the fact that even an infinitesimally small perturbation can make the real part of the contact conductivity negative on some parts of the electrodes. These differentiability results, in fact, hold without more regularity than is assumed in Section~\ref{sec:forward_problem}. Let us start with the derivative with respect to the internal conductivity: The Fr\'echet derivative of $U$ at $\sigma \in L^\infty_+(\Omega)$ in the direction $\eta \in L^\infty(\Omega)$, i.e.,~$D_\sigma U(\sigma; \eta) \in \C_\diamond^M$, is the second part of the unique solution $(D_\sigma u(\sigma; \eta),D_\sigma U(\sigma; \eta)) \in \mathcal{H}^1$ to
\begin{equation}
\label{eq:sderiv0}
B_{\sigma,\zeta}\big((D_\sigma u(\sigma; \eta),D_\sigma U(\sigma; \eta)),(v,V)\big)  \,=  \,  - \int_{\Omega} \eta \nabla u \cdot \nabla v \, {\rm d} x \qquad {\rm for} \ {\rm all} \ (v,V) \in  \mathcal{H}^1,
\end{equation}
where $(u, U) \in \mathcal{H}^1$ is the solution to \eqref{eq:weak}. From the computational standpoint, it is important to note that this derivative can be assembled using the formula (see,~e.g.,~\cite{Hyvonen18})
\begin{equation}
\label{eq:sderiv}
D_\sigma U(\sigma; \eta) \cdot \tilde{I} = - \int_{\Omega} \eta \nabla u \cdot \nabla \tilde{u} \, {\rm d} x,
\end{equation}
 where $(\tilde{u}, \tilde{U}) \in \mathcal{H}^1$ is the solution to \eqref{eq:weak} for an auxiliary current pattern $\tilde{I} \in \C^M_\diamond$. In the following, all solutions for variational equations marked with $\sim$ have $\tilde{I}$ as the underlying electrode current pattern, and the ones without $\sim$ are associated to $I$. The Fr\'echet derivative of $U$ with respect to the contact conductivity allows a similar sampling formula~\cite{Darde21},
\begin{equation}
\label{eq:zderiv}
D_\zeta U(\zeta; \omega) \cdot \tilde{I} = - \int_{\partial \Omega} \omega (U-u) (\tilde{U} - \tilde{u}) \, {\rm d} S
\end{equation}
for $\zeta \in \mathcal{Z}$ and a surface conductivity perturbation $\omega \in L^\infty(E)$.

The map $U: L^\infty_+(\Omega) \times \mathcal{Z} \cap H^1(\partial \Omega) \times \mathcal{B}_c \times \C^M_\diamond \to \C^M_\diamond$ is also Fr\'echet differentiable with respect to its third variable at the origin~\cite{Darde12,Hyvonen17b}. The associated Fr\'echet derivative $D_aU(0; h)$ in the direction $h \in \mathcal{C}^1(E, \R^3)$ is defined by the second part of the unique solution $(D_a u(0; h),D_a U(0; h)) \in \mathcal{H}^1$ to
\begin{equation}
\label{eq:sderiv_a0}
B_{\sigma,\zeta}\big((D_a u(0; h),D_a U(0; h)),(v,V)\big)  \,=  \,  \int_{\partial \Omega} h_\tau \cdot {\rm Grad}(\zeta) \, (U - u) (V  - v) \, {\rm d} S 
\end{equation}
for all $(v,V) \in  \mathcal{H}^1$. Here $h_\tau$ is the tangential component of $h$ and ${\rm Grad}$ is the surface gradient~\cite{Colton98}. This derivative also allows a sampling formula similar to \eqref{eq:sderiv} and \eqref{eq:zderiv} for its efficient numerical evaluation, but it is omitted since it is not utilized in our numerical studies.

Let us finally consider the mixed derivative $D_a D_\sigma U$ that is a main building block in our OED algorithm. To be precise, we consider $D_\sigma D_a U$ without proving that the order of differentiation can be reversed. Moreover, we focus our attention on conductivity perturbations that are supported at a distance from $\partial \Omega$; this choice conveniently follows the material in \cite{Hyvonen14}, but such perturbations are also sufficient for our purposes because the ROI in our numerical experiments always lies in the brain tissue,~i.e.,~in the interior of the computational head model. If we define
$$
L^\infty_\delta(\Omega) = \{ \kappa \in L^\infty(\Omega) \, \big| \, {\rm dist} ({\rm supp} \, \kappa, \partial \Omega) \geq \delta \}, \qquad \delta > 0,
$$
to be the space of admissible conductivity perturbations, the second Fr\'echet derivative $D_\sigma D_a U$ is well defined at $(\sigma, \zeta, 0, I)$ for any  $(\sigma, \zeta,I) \in L^\infty_+(\Omega)\cap \mathcal{C}^{0,1}(\overline{\Omega}) \times \mathcal{Z} \cap H^1(\partial \Omega) \times \C^M_\diamond$. Moreover, it can be assembled via
\begin{align}
  \label{eq:a_sigma_deriv0}
  D_\sigma D_a U(\sigma, \zeta, 0; \eta, h) \cdot \tilde{I} = & \int_{\partial \Omega} h_\tau \cdot {\rm Grad}(\zeta) \big( (D_\sigma U(\sigma; \eta) - D_\sigma u(\sigma; \eta)) (\tilde{U} - \tilde{u}) \nonumber
  \\ & \qquad \quad + (U - u) (D_\sigma \tilde{U}(\sigma; \eta) - D_\sigma \tilde{u}(\sigma; \eta)) \big)  {\rm d} S,
\end{align}
where the perturbation directions satisfy $\eta \in L^\infty_\delta(\Omega)$ and $h \in \mathcal{C}^1(E, \R^3)$. The corresponding formula \cite[(2.21)]{Hyvonen14} for the case of piecewise constant contact conductivity can be deduced from \eqref{eq:a_sigma_deriv0} by noting that $h_\tau \cdot {\rm Grad}(\zeta)$ corresponds to a weighted delta distribution on $\partial E$ for a piecewise constant $\zeta$.

A computational issue with \eqref{eq:a_sigma_deriv0} is that to evaluate $D_\sigma D_a U(\sigma, \zeta, 0; \eta, h)$ for all $\eta$ corresponding to a discretization of the internal conductivity, the variational problem \eqref{eq:sderiv0} must first be solved for each basis current and for as many right-hand sides as there are degrees of freedom in the parametrization of $L^\infty_\delta(\Omega)$. On the other hand, \eqref{eq:a_sigma_deriv0} can be used to conveniently handle perturbations of the electrode positions as one simply needs to plug in \eqref{eq:a_sigma_deriv0} as many different $h$ as there are admissible directions to which the electrodes can move. As the electrode positions are parametrized by their polar and azimuthal angles in our numerical experiments, there are only $2M$ ways to move the electrodes, which is significantly less than the number of basis functions in the parametrization of the conductivity perturbations. To conclude, it would thus be computationally advantageous to replace \eqref{eq:a_sigma_deriv0} by a formula that would not involve solving \eqref{eq:sderiv0} for all conductivity perturbations but instead \eqref{eq:sderiv_a0} for all admissible electrode perturbation directions.

For this reason, interpret the two terms on the right-hand side of \eqref{eq:a_sigma_deriv0} as right-hand sides for \eqref{eq:sderiv_a0} with the test function pairs $(D_\sigma u(\sigma; \eta), D_\sigma U(\sigma;\eta))$ and $(D_\sigma \tilde{u}(\sigma; \eta), D_\sigma \tilde{U}(\sigma;\eta))$, respectively. By the symmetry of the bilinear form $B_{{\sigma},{\zeta}}: \mathcal{H} \times  \mathcal{H} \to \C$, this leads to
\begin{align*}
D_\sigma D_a U(\sigma, \zeta, 0; \eta, h) \cdot \tilde{I} & = {B_{{\sigma},{\zeta}}\big((D_\sigma \tilde{u}(\sigma;\eta), D_\sigma \tilde{U}(\sigma; \eta)), ({D_a u(0; h)}, {D_a U(0; h)}) \big)} \\ & \quad + {B_{{\sigma},{\zeta}}\big((D_\sigma u(\sigma;\eta), D_\sigma U(\sigma; \eta)), ({D_a \tilde{u}(0; h)},{D_a \tilde{U}(0; h)})\big)}.
\end{align*}
The two terms on the right-hand side can now be interpreted as left-hand sides for \eqref{eq:sderiv0}, which finally yields the computationally more attractive formula
\begin{align}
  \label{eq:a_sigma_deriv}
D_\sigma D_a U(\sigma, \zeta, 0; \eta, h) \cdot \tilde{I} & = - \int_{\Omega}  \eta \big( \nabla D_a u(0;h) \cdot \nabla \tilde{u} + \nabla u \cdot \nabla D_a \tilde{u}(0;h) \big)  {\rm d} x
\end{align}
for building a discretized version of the derivative $D_\sigma D_a U$.

\subsection{Head model and computational implementation}

\begin{figure}[t]
\center{
  {\includegraphics[width=6cm]{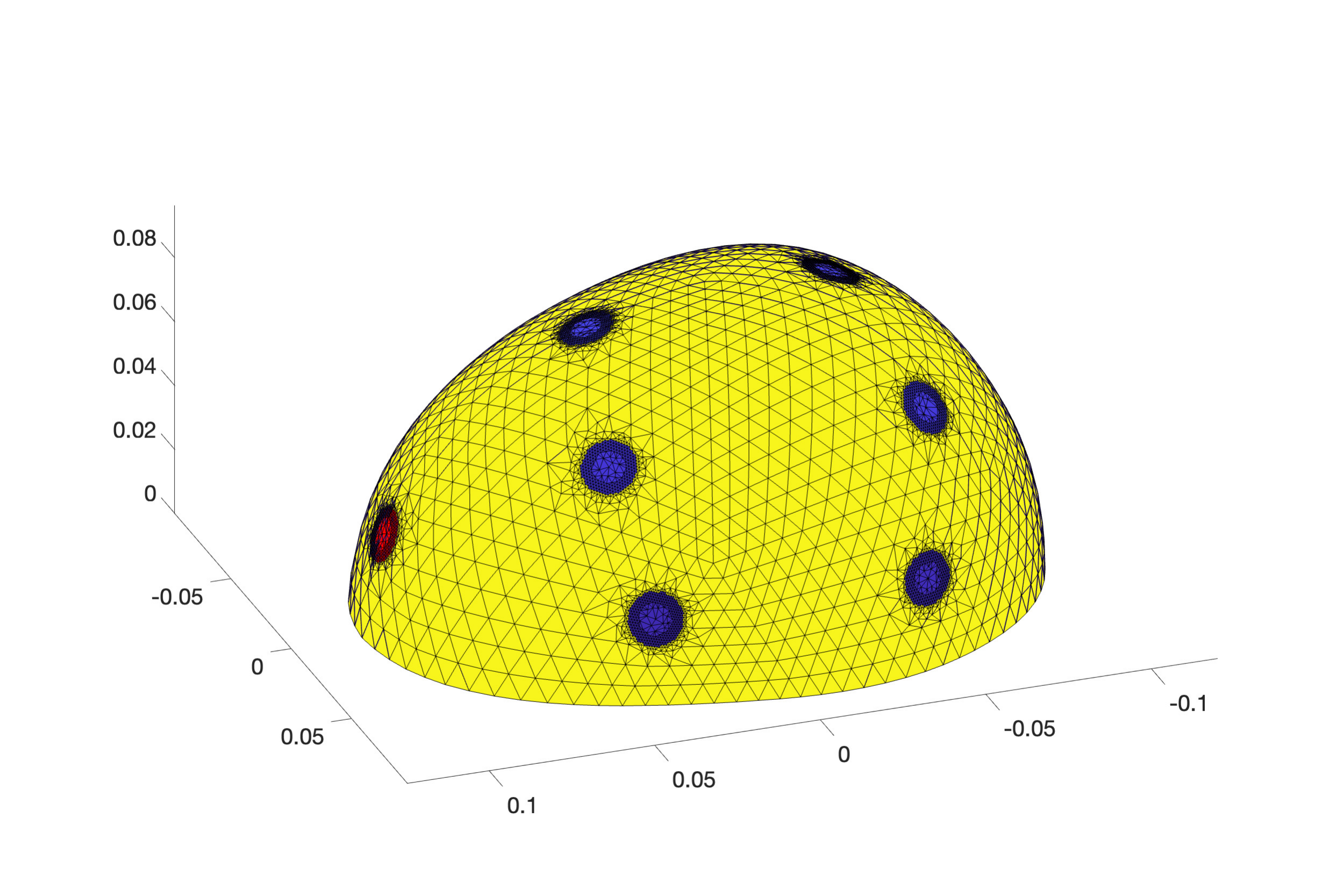}}
  \quad
  {\includegraphics[width=6cm]{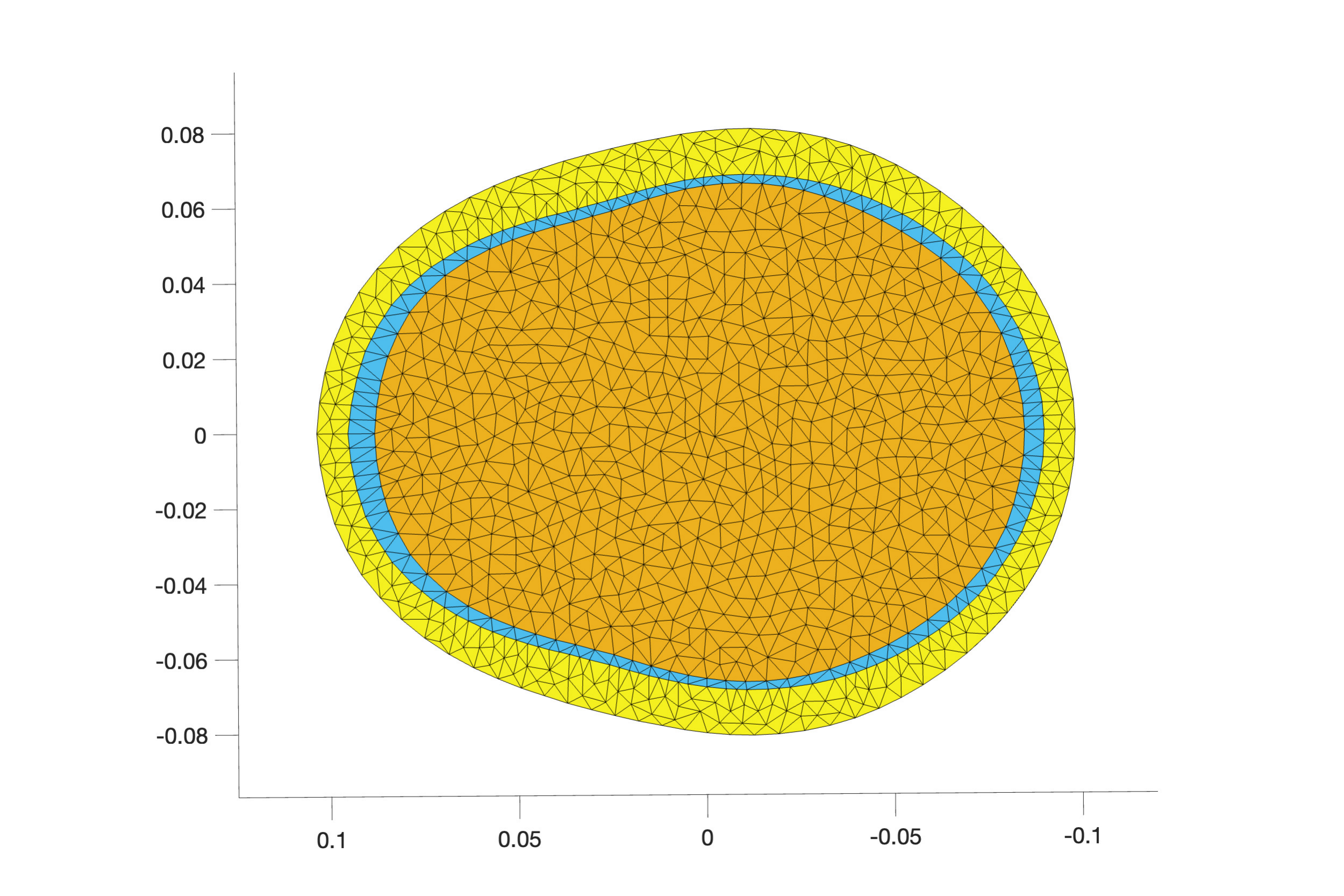}}
  }
\caption{A finite element discretization of the sample head anatomy used in our numerical studies. Left: the surface of the head model with a set of electrodes. Right: the bottom of the computational model showing the skin, skull and brain layers. The unit of length is meter.}
\label{fig:head}
\end{figure}

We adopt the three-layer head model introduced in~\cite{Candiani19, Candiani20}, with the layers corresponding to the scalp, the highly resistive skull and the interior brain. The construction in \cite{Candiani19, Candiani20} models typical variations in the head anatomy over the human population by building a principal component model based on the library of $50$ heads in \cite{Lee16}. However, here we restrict our attention to a single head anatomy, shown in Figure~\ref{fig:head}, without studying the effect of anatomical variations on the optimal electrode configurations. All electrodes shown in Figure~\ref{fig:head} and employed in our numerical studies are discoidal in the sense that they are defined by intersections of the head surface with a circular cylinder of a given radius $R>0$ and having its central axis aligned with the normal of the surface.

Before reviewing how a {\em finite element method} (FEM) can be used to approximate the solution for \eqref{eq:weak} and the derivatives of its second component with the help of \eqref{eq:sderiv}, \eqref{eq:zderiv}, \eqref{eq:sderiv_a0} and \eqref{eq:a_sigma_deriv}, the choice of the perturbation vector field $h$ in \eqref{eq:sderiv_a0} must be considered. We follow the construction in~\cite{Candiani19}. An electrode position is parametrized by the polar $\theta_m \in (0, \pi/2)$ and azimuthal $\phi_m \in [0, 2\pi)$ angles of its center with respect to the midpoint of the bottom face of the head model in Figure~\ref{fig:head}. The tangent vectors for $\partial \Omega$ at the center $x_m$ of $E_m$ obtained by differentiating the parametrization of the head surface (cf.~\cite{Candiani19}) with respect to the polar and azimuthal angles are denoted by $\hat{\theta}_m, \hat{\phi}_m \in \R^3$, respectively.

We extend $\hat{\theta}_m$ as a constant vector field onto the whole electrode $E_m$ and set
$$
\tilde{h}^\theta_m(x) := \big( \hat{\theta}_m(x) \big)_\tau, \qquad x \in E_m,
$$
to be its tangential component. To ensure that all points on $E_m$ move as much in the tangential direction so that the area of the electrode is approximately conserved, we apply a scaling
$$
h^\theta_m(x) := \frac{| \tilde{h}^\theta_m(x_m) |}{| \tilde{h}^\theta_m(x) |} \, \tilde{h}^\theta_m(x), \qquad x \in E_m,
$$
with $| \, \cdot \, |$ denoting the Euclidean norm, and finally extend $h^\theta_m$ as zero on to the other electrodes. Although the perturbation of $E_m$ by such a vector field $h^\theta_m$ in the spirit of \eqref{eq:perturbed} does not exactly maintain the shape of the electrode $E_m$, it seems feasible that the variations in its size are small compared to the moved distance for a regular spherical surface and small electrodes depicted in Figure~\ref{fig:head}. A perturbation vector field $h^\phi_m$ corresponding to the azimuthal angle of the center for $E_m$ is constructed in an equivalent manner. The derivatives of solutions to \eqref{eq:weak} with respect to the polar or the azimuthal angle of an electrode $E_m$  are approximated by solving \eqref{eq:sderiv_a0} with $h = h^\theta_m$ or $h = h^\phi_m$, respectively, for $(D_a u(0; h),D_a U(0; h))$.

In order to parametrize the dependence of the solution to \eqref{eq:weak} on the contact conductivity, we assume the contact takes a common fixed form on each electrode. To be more precise, following~\cite[Section~4.2]{Candiani19}, we define the restriction of $\zeta$ to $E_m$, $m=1, \dots, M$, as  
\begin{equation}
\label{eq:zeta_param}
\zeta|_{E_m}(r_m, \omega_m) =  \zeta_m \hat{\zeta} (r_m, \omega_m) , \qquad r_m \in [0, R), \ \ \omega_m \in [0, 2 \pi),
\end{equation}
in the polar coordinates $(r_m, \omega_m)$ induced by orthogonally projecting the electrode onto a cross-section of the circular cylinder defining it. As in \cite{Candiani19}, the shape of an individual contact is defined by the smooth function
$$
\hat{\zeta} (r, \omega) = \exp\left( \tau - \frac{\tau R^2}{R^2 - r^2} \right)
$$
that vanishes to infinite order at $|r| = R$. In our numerical experiments, the value $\tau = 0.4$ is used for the shape parameter. The free parameter $\zeta_m \in \C$, with ${\rm Re}\, \zeta_m > 0$, in \eqref{eq:zeta_param} defines the net conductance on $E_m$. In what follows, we abuse the notation by identifying the contact conductance function with the vector $\zeta = [\zeta_1, \dots, \zeta_M]^\top \in \C_+^M$, where $\C_+$ denotes the open right half of the complex plain.

In the rest of this work, we consider a discretized version of \eqref{eq:weak} by parametrizing the conductivity $\sigma$ as
\begin{equation}
\label{eq:discr_sigma}
\sigma = \sum_{j=1}^N \sigma_j \varphi_j \, ,
\end{equation}
where $\sigma_j \in \C_+$ and $\varphi_j \in H^1(\Omega)$, $j=1, \dots, N$,  are the piecewise linear basis functions corresponding to a {\em finite element} (FE) mesh of $\Omega$.\footnote{In fact, the basis functions $\varphi_j$ and their number $N$ change when electrodes move on $\partial \Omega$ due to remeshing.} We refer to \cite{Candiani19, Candiani20} for more information on the construction of tetrahedral FE meshes such as shown in Figure~\ref{fig:head}. We denote by $\sigma$ both the vector $\sigma \in \C^N_+$ and the corresponding conductivity defined by \eqref{eq:discr_sigma}; such an identification between linear combinations of the FE basis functions and vectors of  $\C^N$ is employed throughout this text. By the solution $(u,U)$ of~\eqref{eq:weak} we mean from here on the approximate FE solution in the subspace
$$
{\rm span}\{\varphi_j\}_{j=1}^N \oplus \C^M_\diamond \subset H^1(\Omega)\oplus \C^M_{\diamond};
$$
consult,~e.g.,~\cite{Vauhkonen97, Vauhkonen99} for the implementation details.

Summarizing the above developments, we (re)define a finite-dimensional forward map 
\begin{equation}
  \label{eq:inv_forv}
  \mathcal{U}:
  \left\{
  \begin{array}{l}
    (\sigma, \zeta, \theta, \phi) \to \mathcal{U}(\sigma, \zeta, \theta, \phi), \\[1mm]
    \C_+^{N} \times \C^M_+ \times (0,\pi/2)^M \times [0, 2\pi)^M \to \C^{M(M-1)}
  \end{array}
  \right.
\end{equation}
that is a discretized version of  \eqref{eq:forward_map} with the dependence on the electrode angles explicitly accounted for. The Jacobian of $J_\sigma = J_\sigma(\sigma, \zeta, \theta, \phi) \in \C^{M(M-1) \times N}$ with respect to $\sigma \in \C_+^{N}$ can be computed by letting $\eta$ run through the FE basis functions in \eqref{eq:sderiv} for all $I$ and $\tilde{I}$ in the employed current basis. Similarly, the Jacobian $J_\zeta = J_\zeta(\sigma, \zeta, \theta, \phi) \in \C^{M(M-1) \times M}$ with respect to $\zeta \in \C^M_+$ can be assembled by choosing $\omega$ in \eqref{eq:zderiv} to be of the form \eqref{eq:zeta_param} with $\zeta$ going through the standard basis vectors of $\C^M$. Finally, the elementwise derivatives
\begin{equation}
  \label{eq:second_deriv}
\frac{\partial (J_\sigma)_{ij}}{\partial \theta_m}(\sigma, \zeta, \theta, \phi) \qquad \text{and} \qquad  \frac{\partial (J_\sigma)_{ij}}{\partial \phi_m}(\sigma, \zeta, \theta, \phi)
\end{equation}
for $i=1, \dots, M(M-1)$, $j = 1,\dots, N$ and $m=1, \dots, M$, needed in the OED algorithm for the linearized model, are approximated by setting $h = h^\theta_m$ or $h = h^\phi_m$ on the right-hand side of \eqref{eq:a_sigma_deriv}, choosing $\eta = \varphi_j$, setting the current pattern $I$ to be the one corresponding to the index $i$ in $\mathcal{U}(\sigma, \zeta, \theta, \phi)$ (cf.~\eqref{eq:forward_map}), and letting $\tilde{I}$ run through the whole current basis of $\C^M_\diamond$.

Let us then assume the knowledge of background conductivity levels $\sigma_0 \in \C_+^N$ and $\zeta_0 \in \C_+^M$. In the stroke monitoring framework, the former would be based on information on the physiology of the patient, whereas the latter would correspond to empirical information on the expected quality of electrode contacts. We resort to linearizing the map $\mathcal{U}: \C_+^{N} \times \C^M_+ \times (0,\pi/2)^M \times [0, 2\pi)^M \to \C^{M(M-1)}$ with respect to either its first or both its first and second variables. Assuming a noiseless measurement $\mathcal{V} \in \C^{M(M-1)}$ and performing the linearization with respect to both $\sigma$ and $\zeta$ leads to the approximate measurement model
  \begin{equation*}
    \mathcal{U}(\sigma_0, \zeta_0) + J_\sigma(\sigma_0, \zeta_0) (\sigma - \sigma_0) + J_\zeta(\sigma_0, \zeta_0) (\zeta - \zeta_0) = \mathcal{V},
  \end{equation*}
  where the dependence of the Jacobians on the design parameters $\theta$ and $\phi$ has been suppressed for brevity. Denoting $\kappa = \sigma - \sigma_0$, $\xi = \zeta - \zeta_0$, $y = \mathcal{V} - \mathcal{U}(\sigma_0, \zeta_0)$,
  \begin{equation}
    \label{eq:J_and_w}
   J = [J_\sigma(\sigma_0, \zeta_0),  J_\zeta(\sigma_0, \zeta_0)] \quad \text{and} \quad  w =  \begin{bmatrix}
    \kappa \\
    \xi
    \end{bmatrix},
\end{equation}
  we arrive at the model
  \begin{equation}
        \label{eq:linearization}
  J w = y
   \end{equation}
  that is further analyzed from Bayesian standpoint in the next section.

  When considering the actual problem of OED, we assume that $\zeta_0$ is a good enough estimate for the contacts and set $J = J_\sigma$ and $w = \kappa$ in \eqref{eq:linearization}. That is, the complete linearized model \eqref{eq:linearization} is only considered when computing reconstructions with the TV prior by the algorithm based on the lagged diffusivity iteration. The reason for this choice is two-fold: we do not expect the dependence of the model on $\zeta$ to severely affect the optimal designs aiming at reconstructing $\sigma$ and we want to avoid introducing extra technical material in the form of another second derivative $D_\zeta D_a U$ in addition to $D_\sigma D_a U$. In the following, we denote the design variables by $d = (\theta, \phi)$ and stress that $J = J(d)$.

\section{Bayesian inversion}
\label{sec:Bayesian}

In this section, we first describe a fully Gaussian probabilistic version of the linearized model \eqref{eq:linearization}, with the dependence on the contact conductivity suppressed. Subsequently, a smoothened TV prior and the lagged diffusivity iteration are considered. For the rest of this text, we assume for simplicity that all variables are real-valued.

\subsection{A fully Gaussian model for OED}
\label{sec:Gaussian}

 Let us consider a probabilistic version of \eqref{eq:linearization} with the dependence on the contact conductivity ignored:
\begin{equation}
  \label{eq:meas_model}
Y = J_\sigma W + E,
\end{equation}
where the measurement data $Y \in \R^{M(M-1)}$ and the unknown {\em internal} conductivity perturbation $W \in \R^{N}$ are interpreted as independent random variables. The additive noise process $E$ is assumed to be a zero-mean Gaussian with a symmetric positive-definite covariance matrix $\Gamma_{{\rm noise}} \in \R^{M(M-1) \times M(M-1)}$.

Let us assume that {\em a priori} $W \sim \mathcal{N}(0, \Gamma_{\rm prior})$, where the prior covariance $\Gamma_{\rm prior} \in \R^{N \times N}$ is also symmetric and positive definite. In such a case, it is well known that the posterior density $\pi(w \, | \, y)$ is a Gaussian as well (see,~e.g.,~\cite{Kaipio06}), with the covariance
  \begin{align}
  \label{eq:Gaussian_posterior}
    \Gamma_{\rm post} =  \Gamma_{\rm prior} - \Gamma_{\rm prior} J_\sigma^\top\big(J_\sigma \Gamma_{\rm prior} J_\sigma^\top + \Gamma_{\rm noise} \big)^{-1} J_\sigma \Gamma_{\rm prior}
   \end{align}
and the mean 
\begin{align}
\label{eq:Gaussian_posterior_mean}
  \widehat{w} =\Gamma_{\rm prior} J_\sigma^\top\big(J_\sigma \Gamma_{\rm prior} J_\sigma^\top + \Gamma_{{\rm noise}} \big)^{-1} y.
  \end{align}
In particular, the posterior covariance does not depend on the measured data $y$, and thus an optimal design aiming to make the posterior as localized as possible can be sought offline,~i.e.,~prior to taking any measurements.

In our numerical experiments on OED for the head model, a weighted trace of the posterior covariance $\Gamma_{\rm post} = \Gamma_{\rm post}(d)$ is minimized with respect to the design variable $d$. Note that the dependence of $\Gamma_{\rm post}$ on $d$ is inherited from that of $J_\sigma = J_\sigma(d)$. The prior covariance in \eqref{eq:Gaussian_posterior} originates either from {\em a priori} assumption on the conductivity perturbation, or it is a byproduct of a reconstruction algorithm based on sequential linearizations and the lagged diffusivity iteration, as described in the following section.

\subsection{Total variation prior and lagged diffusivity}
\label{sec:TV}

Let us then assume that the prior density is of the form
\begin{equation}
  \label{eq:TV}
\pi(w) = \exp(- \gamma \Phi(\kappa)),
\end{equation}
where 
\begin{equation}
\label{eq:TVexp}
\Phi(\kappa) = \int_{\Omega} \upsilon(x) \varphi \big(|\nabla \kappa | \big) \, {\rm d} x, \qquad \text{with} \quad \varphi(t) = \sqrt{t^2 + T^2} \approx | t |,
\end{equation}
accompanied by the information that the interior conductivity perturbation $\kappa$ vanishes on $\partial \Omega$. As in \cite{Candiani21}, $\upsilon: \Omega \to \R_+$ is the reciprocal of a smooth cut-off function,
\begin{equation}
    \label{eq:upsilon}
\upsilon(x) = \Big( \frac{1}{2}\big( 1 + \tanh( c_\upsilon  ({\rm dist}(x, \partial \Omega) - b_\upsilon) ) \big) \Big)^{-1}, \qquad x \in \Omega,
\end{equation}
where $c_\upsilon > 0$ and $b_\upsilon > 0$ are chosen so that the value of $\upsilon$ is large in the skin layer and decreases to almost one within the skull layer. The role of $\upsilon$ in the prior is to encode the assumption that a stroke mainly causes conductivity changes in the brain tissue, but in the computational experiments it also mitigates the tendency of electrodes to approach the bottom face of the head model because $\upsilon$ also makes changes in the conductivity improbable in the lower parts of the brain. Note that the vector of perturbations in the peak values of the contact conductivity, $\xi$, is assumed to have an uninformative prior, which explains why $\xi$ does not appear on the right-hand side of \eqref{eq:TV}. In particular, the dependence of the linearized model \eqref{eq:linearization} on $\xi$ is not ignored in this section. We assume that there are $N' < N$ interior nodes in the employed FE mesh and identify $\kappa$ in the rest of this section with an element of $\R^{N'}$, instead of one in $\R^N$. 

Let us continue to assume a zero-mean Gaussian noise model that is independent of $W$ and follow the treatment in \cite{Harhanen15}. The quest for a MAP estimate for the linearized model \eqref{eq:linearization} leads now to considering  a nonquadratic Tikhonov functional
\begin{equation}
\label{Tikhonovk}
 \frac{1}{2}\big(y - J w )^{\rm \top} \Gamma_{\rm noise}^{-1} \big(y - J w \big) +  \gamma \, \Phi(\kappa).
\end{equation}
The necessary condition for a minimizer of \eqref{Tikhonovk} is (cf.~\eqref{eq:J_and_w})
\begin{equation}
\label{eq:necessary}
J^{\rm \top} \Gamma_{\rm noise}^{-1} J  w  
+ \gamma
\begin{bmatrix}
(\nabla_\kappa \Phi)(\kappa) \\[1mm]
\mathrm{0}
\end{bmatrix}
 =  J^{\rm \top} \Gamma_{\rm noise}^{-1}y,
\end{equation}
where $\mathrm{0} \in \R^M$. Utilizing the identification of vectors in $\R^{N'}$ with elements of $H^1_0(\Omega)$ via the employed FE basis, it is straightforward to check that
\begin{equation}
  \label{eq:kappa_deriv}
(\nabla_\kappa \Phi)(\kappa) = \Theta(\kappa) \kappa,
\end{equation}
where the multiplier matrix $\Theta(\kappa)$ is given elementwise as
$$
\Theta_{i,j}(\kappa) := \int_{\Omega} \frac{\upsilon(x)}{\sqrt{|\nabla_x \kappa (x)|^2 + T^2}} \nabla \varphi_i(x) \cdot \nabla \varphi_j(x) \, {\rm d} x, \qquad i,j=1,\dots ,N'.
$$
As reasoned in \cite{Harhanen15}, the invertibility of $\Theta(\kappa)$ for any $\kappa \in \R^{N'}$ is guaranteed via an interpretation as a FE system matrix originating from a discretization of an elliptic differential operator combined with a homogeneous Dirichlet boundary condition.

Plugging the presentation \eqref{eq:kappa_deriv} in \eqref{eq:necessary} and expanding gives
\begin{equation}
\label{eq:necessary2}
\begin{bmatrix}
J_\sigma^{\top} \Gamma_{\rm noise}^{-1} J_\sigma + \gamma \, \Theta(\kappa) & J_\sigma^{\top} \Gamma_{\rm noise}^{-1} J_\zeta  \\[1mm]
J_\zeta^{\top} \Gamma_{\rm noise}^{-1} J_\sigma & J_\zeta^{\top} \Gamma_{\rm noise}^{-1} J_\zeta \!\!
\end{bmatrix}
 w
 =  \begin{bmatrix}
   J_\sigma^{\top} \\[1mm] J_\zeta^\top
   \end{bmatrix}
   \Gamma_{\rm noise}^{-1}\, y.
\end{equation}
Denoting $B_\sigma = \Gamma_{\rm noise}^{-1/2} J_\sigma$, $B_\zeta = \Gamma_{\rm noise}^{-1/2} J_\zeta$ and solving the second row of \eqref{eq:necessary2} for $\xi$ gives
\begin{equation}
\label{eq:xi_update}
  \xi = \big(B_\zeta^{\top} B_\zeta \big)^{-1} B_\zeta^{\top} \big( \Gamma_{\rm noise}^{-1/2} y  - B_\sigma \kappa \big).
  \end{equation}
Substituting this for $\xi$ on the first row of \eqref{eq:necessary2} results in the formula
\begin{equation}
  \label{eq:normal_form}
\big( A^{\rm T}A + \gamma \Theta(\kappa) \big) \kappa = A^{\rm T}b, 
\end{equation}
where
$$
A =  Q B_\sigma \qquad {\rm and} \qquad b  =  Q \Gamma_{\rm noise}^{-1/2} y, 
$$
with $Q = \mathrm{I} - B_\zeta (B_\zeta^{\top}  B_\zeta)^{-1} B_\zeta^{\top}$ being the orthogonal projection onto the orthogonal complement of the range of $B_\zeta$ in $\R^{M(M-1)}$. Here and in what follows, an identity matrix of an appropriate size is denoted by $\mathrm{I}$.

The idea is now to solve \eqref{eq:normal_form} by means of the lagged diffusivity iteration \cite{chan1999convergence,dobson1997convergence,Vogel96}, which can be interpreted to produce a Gaussian approximation for the posterior of $\kappa$ as its byproduct~\cite{Bardsley18,Calvetti08,Helin22}. Start from an initial guess $\kappa^{0} = 0$ for the solution of \eqref{eq:normal_form}. Assuming to have $\kappa^{(j)}$ in hand, $\Theta(\kappa)$ is replaced in \eqref{eq:normal_form} by $\Theta(\kappa^{(j)})$, which yields a linear equation for $\kappa$. The solution of this equation is the mean of the Gaussian posterior defined by the model for $A \kappa = b$, assuming a zero-mean Gaussian prior for $\kappa$ with the inverse covariance matrix $\gamma \Theta(\kappa^{(j)})$ and an additive white noise model~\cite{Kaipio06}. Using the Woodbury matrix identity (cf.~\eqref{eq:Gaussian_posterior}--\eqref{eq:Gaussian_posterior_mean}), the mean and the covariance of this posterior allow the representations
\begin{equation}
  \label{eq:LD_mean}
\kappa^{(j+1)} =  \Theta(\kappa^{(j)})^{-1} A^\top \big(\gamma \mathrm{I} + A \Theta(\kappa^{(j)})^{-1} A^\top \big)^{-1} b
\end{equation}
and
\begin{equation}
  \label{eq:LD_covariance}
\Gamma^{(j+1)} = \gamma^{-1} \Big( \Theta(\kappa^{(j)})^{-1}   -
\Theta(\kappa^{(j)})^{-1} A^\top \big(\gamma \mathrm{I} + A \Theta(\kappa^{(j)})^{-1} A^\top \big)^{-1} A \Theta(\kappa^{(j)})^{-1} \Big),
\end{equation}
respectively. 

The iteration \eqref{eq:LD_mean} is terminated after a predefined stopping criterion is satisfied. The final posterior mean and covariance, say, $\kappa^{(J)}$ and $\Gamma^{(J)}$, then define a Gaussian density that can be used in choosing new positions for the electrodes as explained in Section~\ref{sec:oed}. Note that the iteration defined by \eqref{eq:LD_mean} is not overly expensive computationally because $A$ has only $M(M-1)$ rows and operating with the inverse of the sparse matrix $\Theta(\kappa^{(j)})$ that originates from a FE discretization can be performed efficiently. Moreover, the covariance matrix in \eqref{eq:LD_covariance} only needs to be computed at the final iteration, since \eqref{eq:LD_mean} is independent of \eqref{eq:LD_covariance}. After computing a MAP estimate for the posterior by the above scheme, the peak values of the contact conductivity perturbation could be updated according \eqref{eq:xi_update} with $\kappa = \kappa^{(J)}$.

\begin{remark}
  \label{rmrk:seq_lin}
  In the numerical experiments of Section~\ref{sec:numerics}, we compute the TV-based reconstructions by combining the lagged diffusivity iteration with sequential linearizations of the forward model. That is, once a conductivity perturbation pair, say, $w_k = (\kappa_k, \xi_k)$, has been computed by the above iteration, the forward model is re-linearized at the corresponding new base point $(\sigma_k, \zeta_k) = (\sigma_0 + \kappa_k, \zeta_0 + \xi_k)$ and the lagged diffusivity iteration is once again applied to the nonquadratic Tikhonov functional corresponding to the negative log-likelihood for the TV prior. The unknown $(\kappa_{k+1}, \xi_{k+1})$ that is sought for in the new round of iteration still corresponds to a deviation from the background level $(\sigma_0, \zeta_0)$, which means that one must redefine
  $$
  y = \mathcal{V} - \mathcal{U}(\sigma_k, \zeta_k) + J_\sigma(\sigma_k, \zeta_k) \kappa_k + J_\zeta(\sigma_k, \zeta_k) \xi_k + e
  $$
  in the analysis of this section, with $e$ denoting the considered realization of the additive measurement noise process $E$. The reason for this choice is that it is natural to assume the TV prior for the change in the internal conductivity compared to the original background level $\sigma_0$, not compared to the latest reconstruction $\sigma_k$. Once the nested iterations corresponding to the sequential linearizations and the lagged diffusivity argument have been terminated, the final approximate Gaussian posterior covariance for the perturbation from the background conductivity $\sigma_0$ can be used for choosing new electrode positions as explained in the next section. See \cite{Candiani21,Harhanen15} for previous applications of such a reconstruction algorithm to EIT, with the difference that utilizing the Woodbury matrix formula for \eqref{eq:LD_mean} is in \cite{Candiani21,Harhanen15} replaced by a certain priorconditioned LSQR iteration and the Morozov discrepancy principle.
  \end{remark}

\section{Optimal experimental design}
\label{sec:oed}
In this section, we briefly describe the concept of A-optimality in our setting and consider numerically solving the associated optimization problems. We refer to~\cite{alexanderian2016bayesian,chaloner1995bayesian} for more information on Bayesian OED and to \cite{Nocedal06} for material on numerical optimization.

\subsection{A-optimality}
\label{ssec:A_optimality}
 Assume that the forward model has been linearized around a given internal background conductivity and that the prior information on the conductivity perturbation is given in the form of a zero-mean Gaussian prior with a covariance matrix $\Gamma_{\rm prior}$. In other words, we assume the setting of Section~\ref{sec:Gaussian}. The prior information may originate from one of two sources:
\begin{enumerate}
\item A Gaussian prior is directly assumed for the linearized model, and the computations are performed offline in a nonadaptive manner.
\item The algorithm based on sequential linearizations, a TV prior and the lagged diffusivity iteration, described in Section~\ref{sec:TV}, is first applied to measurements on an initial electrode configuration. After the algorithm has been terminated, the covariance of the final Gaussian approximation for the internal conductivity perturbation is dubbed $\Gamma_{\rm prior}$.
\end{enumerate}

In both cases, the covariance $\Gamma_{\rm post}$ of the Gaussian posterior for the linearized model is given by \eqref{eq:Gaussian_posterior}. The A-optimality criterion corresponds to minimizing the trace of a weighted version of this covariance:
\begin{equation}
\label{eq:Aoptimal}
d_{\rm A} =  \underset{d}{\argmin} \, \, \Psi_{\rm A}(d), \quad \  \text{with} \ \Psi_{\rm A}(d) := {\rm tr}  \big(\mathcal{A} \Gamma_{\rm post} (d) \mathcal{A}^\top\big) =  {\rm tr}\big(\mathcal{A}^\top \! \mathcal{A} \Gamma_{\rm post} (d) \big),
\end{equation}
where the last equality follows from a cyclic property of the matrix trace. In other words, an A-optimal design $d_{\rm A}$ minimizes the expected squared distance of the unknown from the posterior mean in the seminorm defined by the positive semidefinite matrix $\mathcal{A}^\top \! \mathcal{A}$; see,~e.g.,~\cite{Burger21}. In our setting, the weighing provided by $\mathcal{A}$ for the diagonal elements of $\Gamma_{\rm post}$ is essential as they correspond to the nodal values of FE basis functions with supports of different sizes. To be able to measure distances in the squared $L^2(\Omega)$ norm, the proper choice for $\mathcal{A}^\top \!\mathcal{A}$ is thus the mass matrix corresponding to the used FE basis; see, e.g., \cite{Hannukainen20}.  Moreover, if only a certain ROI inside the head model is considered, as is often the case in our numerical tests, the FE basis functions corresponding to the nodes outside the ROI are replaced by zero functions when forming the mass matrix.

When visualizing the evolution of the A-optimality target during its minimization, instead $\Psi_{\rm A}$ itself, the value of $\widetilde{\Psi}_{\rm A} = \Psi_{\rm A}^{1/2}$ is monitored. With the above choice for the weight $\mathcal{A}$, the value of $\widetilde{\Psi}_{\rm A}$ is the square root of the expected squared $L^2$ error over the ROI under the linearized measurement model.

\subsection{Minimization of the A-optimality target}
The minimization of the A-optimality target $\Psi_{\rm A}$ is performed with gradient descent accompanied by an inexact line search based on Armijo's rule; see,~e.g.,~\cite{Nocedal06}. The complete minimization algorithm is given as a combination of Algorithms~\ref{alg:GD} and \ref{alg:armijo_alg}.

Observe that in \eqref{eq:Gaussian_posterior}, only the Jacobian matrix $J_\sigma = J_\sigma(d)$ depends on the design parameter vector $d$ that determines the positions of the electrodes. The elementwise derivatives of  $J_\sigma(d)$ can be approximated using \eqref{eq:second_deriv}. Hence, the derivatives of the optimization target $\Psi_{\rm A}(d)$ needed in Algorithm~\ref{alg:GD} can also be computed straightforwardly, yet tediously, by applying \eqref{eq:second_deriv} and well-known matrix differentiation formulas to \eqref{eq:Gaussian_posterior} and \eqref{eq:Aoptimal}. However, we do not stress this matter any further here, but simply assume that we have access to an approximation of the gradient of $\Psi_{\rm A}$ when running Algorithms~\ref{alg:GD} and \ref{alg:armijo_alg}.

\begin{algorithm}[t]
\caption{(Gradient descent)} 
\label{alg:GD}
\begin{algorithmic} 
\STATE Choose a tolerance $\epsilon > 0$, an initial guess $d_0$, the maximum number of iterations $N_{\rm GD}$, and a step size parameter $\lambda > 0$. Initialize $i = 0$. 
\WHILE {$|\nabla \Psi_{\rm A}(d_i)| > \epsilon$ and $i < N_{\rm GD}$}
\STATE \hspace{2mm} $\rhd$ Compute the search direction $q = \tilde{q}/|\tilde{q}|$, with $\tilde{q} = -\nabla \Psi_{\rm A}(d_i)$.
\STATE \hspace{2mm} $\rhd$ Select the step size $\bar{\lambda}$ based on Algorithm~\ref{alg:armijo_alg} with an initial guess $\lambda$.
\STATE \hspace{2mm} $\rhd$ Set $d_{i+1} = d_i + \bar{\lambda} q$.
   \STATE \hspace{2mm} $\rhd$ Update $i = i+1$.    
\ENDWHILE
\STATE {\bf return} $d_i$
\end{algorithmic}
\end{algorithm}

\begin{algorithm}[t]
\caption{(Inexact line search with Armijo's rule)}
\label{alg:armijo_alg}
\begin{algorithmic}
\STATE Choose the maximum number of iterations $N_{\mathrm{Armijo}}$ and scaling parameters $\alpha, \beta \in (0,1)$. Initialize a counter as $l = 0$. The initial step size $\lambda > 0$, the considered base point $d_i$, and the search direction $q$ are given as inputs.
\WHILE{$l < N_{\mathrm{Armijo}}$}
    \IF{$\Psi_{\rm A}(d_i + \lambda q) - \Psi_{\rm A}(d_i) < \lambda \alpha \nabla \Psi_{\rm A}(d_i)^\top q$}
      \STATE \hspace{2mm} $\rhd$ Return $\bar{\lambda} = \lambda$.
    \ELSE
        \STATE \hspace{2mm} $\rhd$ Set $\lambda = \beta \lambda$.
    \ENDIF
    \STATE \hspace{2mm} $\rhd$ Set $l = l + 1$.
\ENDWHILE
\RETURN $\bar{\lambda}  = \lambda$
\end{algorithmic}
\end{algorithm}

\begin{remark}
When the polar angle $\theta_m$ of the center of the electrode $E_m$ is small,~i.e.,~the electrode lies close to the `north pole' of the head model, the derivative with respect to the azimuthal angle $\phi_m$ of the center of $E_m$ in Algorithms~\ref{alg:GD} and \ref{alg:armijo_alg} is not used to update the azimuthal angle itself but the projected position of the electrode center along a straight line on the tangent plane touching the head model at the north pole. The reason for this is to prevent electrodes with small polar angles from rotating around the north pole during the optimization process.
\end{remark}

\section{Numerical experiments}
\label{sec:numerics}
The computations are performed on the three-layer head model presented in Figure~\ref{fig:head}. The employed FE meshes, used for discretizing both the conductivity and the first component of a solution to \eqref{eq:weak}, have approximately $N=15\,000$ nodes and $70\,000$ tetrahedra with appropriate refinements at the electrodes. The precise characteristics of the mesh depend on the positions of the electrodes on the surface of the head. However, when the positions of the electrodes are altered, the mesh only changes in the skin and skull layers,~i.e.,~the mesh for the brain is the same in all computations. The needed FE meshes are formed using an updated version of the workflow introduced in~\cite{Candiani19,Candiani20}, with emphasis on the robustness of meshing when electrodes are close to each other or near the bottom edge of the head model.

\begin{figure}[t]
\center{
  {\includegraphics[width=6cm]{Figures/evenposl.eps}}
  \quad
      {\includegraphics[width=6cm]{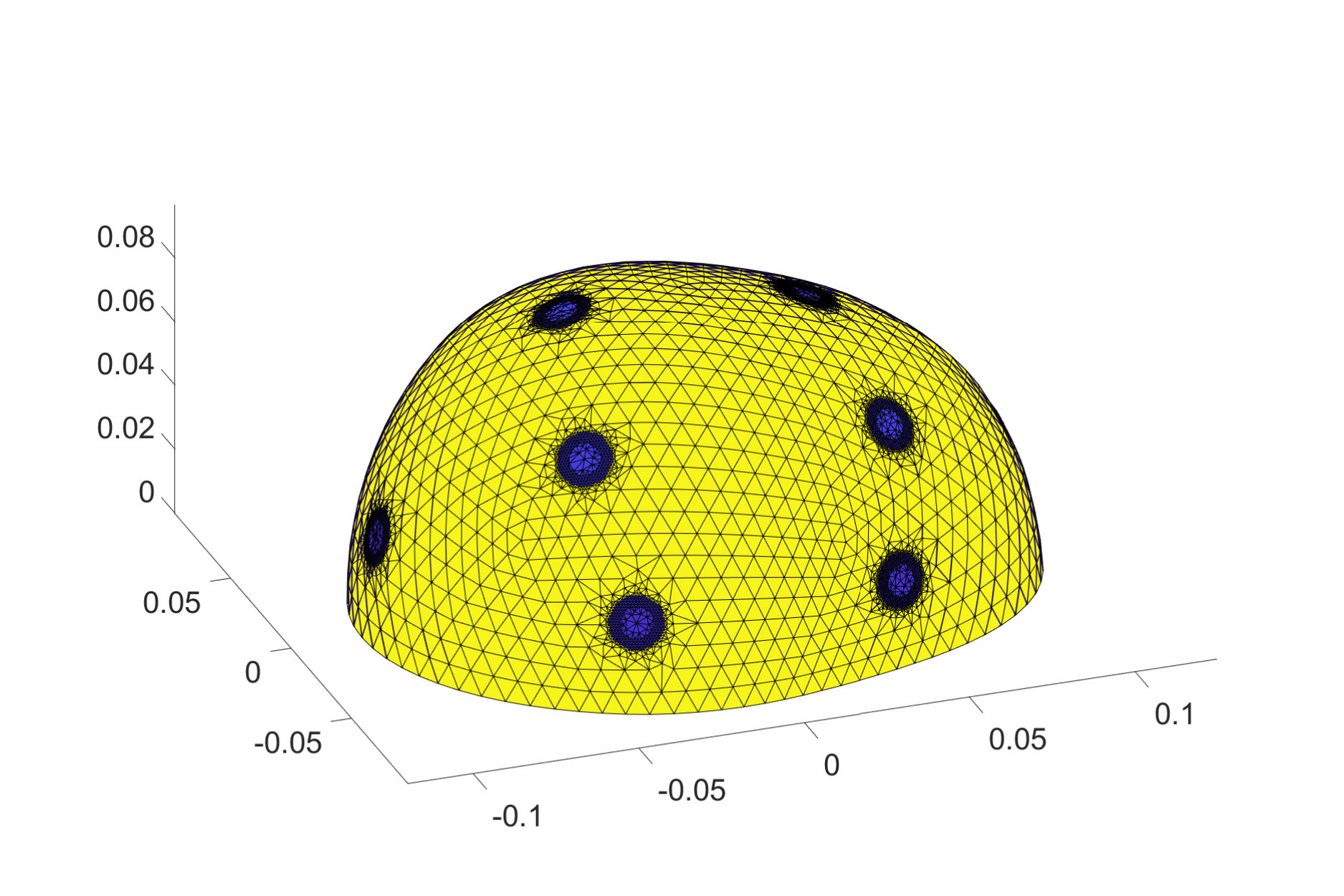}} \linebreak
      {\includegraphics[width=6cm]{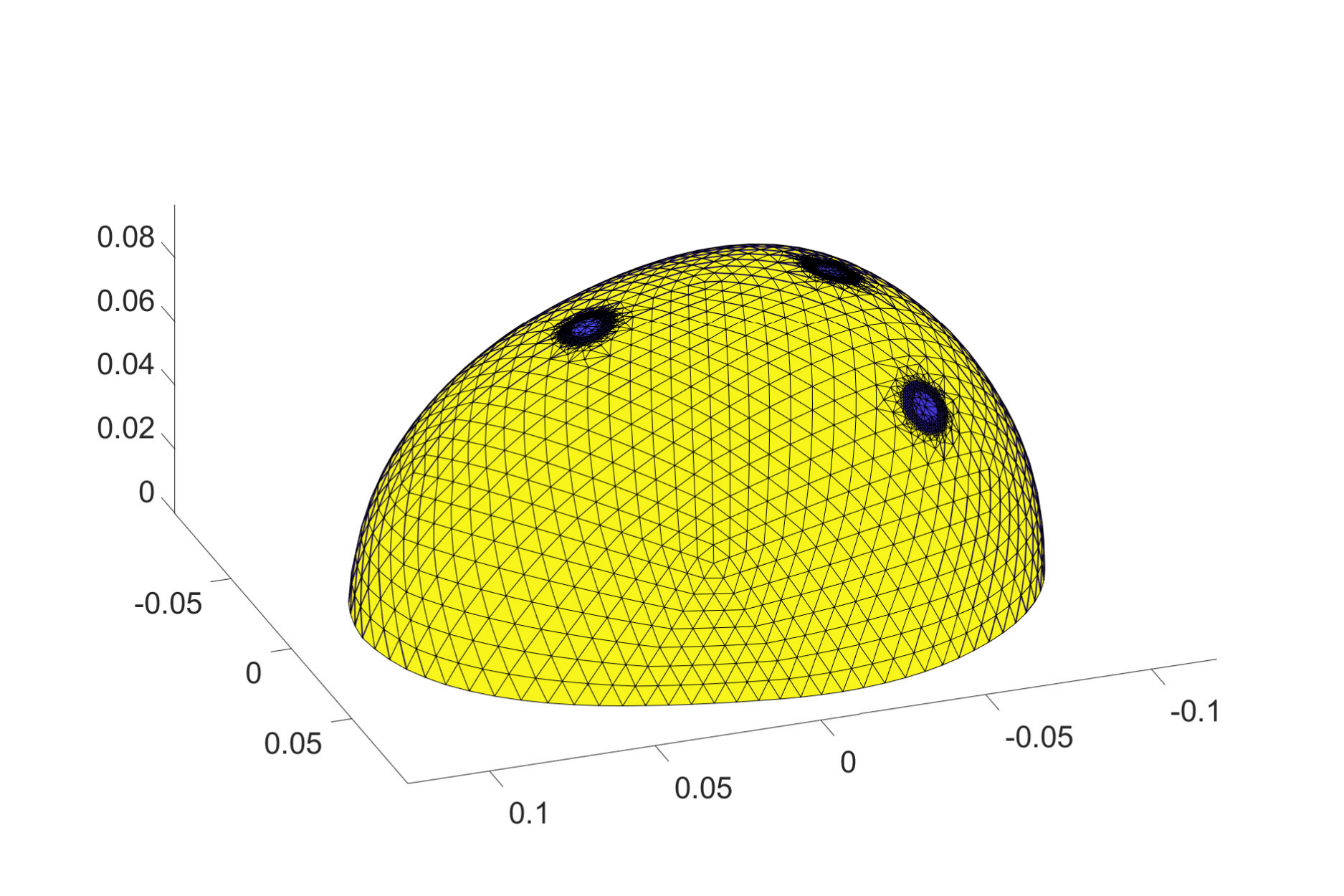}}
  \quad
      {\includegraphics[width=6cm]{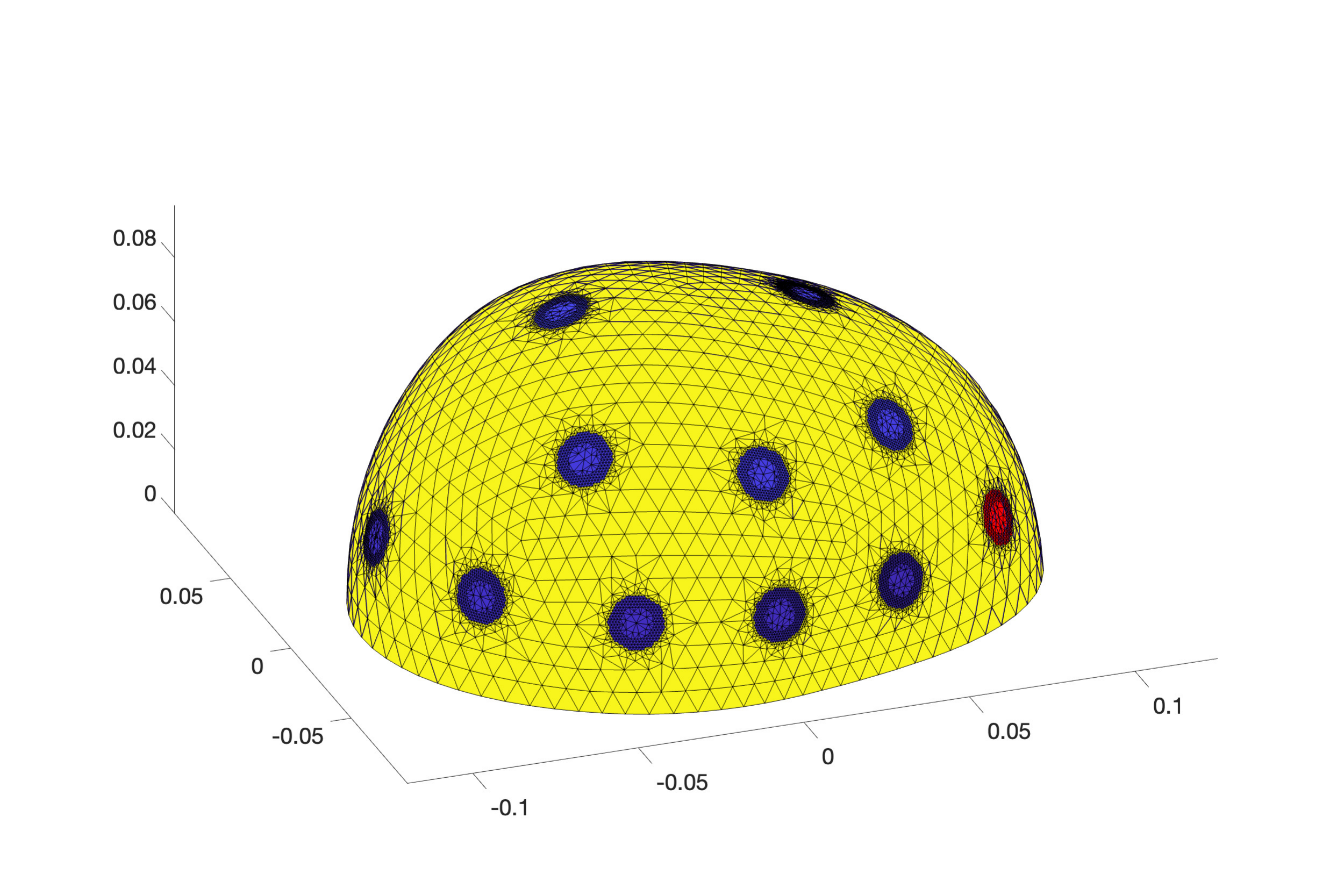}}
  }
\caption{The employed initial configurations with $M=12$ electrodes. The feeding electrode is shown in red. Top: a symmetric initial electrode configuration from two perspectives. Bottom: an initial electrode configuration concentrated close to the ROI or stroke in the right posterior quadrant of the brain from two perspectives.}
\label{fig:init_conf}
\end{figure}

The number of electrodes is chosen to be $M=12$ in all numerical experiments, and their common radius is set to $7.5$\,mm. Denoting the $l$th Cartesian basis vector by ${\rm e}_l$, the employed current basis are $I_m = {\rm e}_1 - {\rm e}_{m+1}$, $m=1, \dots, M-1$, that is, the same electrode always feeds the current into the head and it flows out in turns through the other $M-1$ electrodes. In the numerical tests, we consider the two initial configurations for the electrodes shown in Figure~\ref{fig:init_conf}, with the current feeding electrodes marked in red. Although the effect of the initial position for the feeding electrode is briefly discussed in Section~\ref{sec:gauss_numer}, we refrain from further analysis on the role of the current patterns and settle for referring to \cite{Kaipio04,Kaipio07} for their optimization in the Bayesian framework.

The peak values of the contact conductances are assumed to be known and set to $\zeta_1, \dots, \zeta_M = 10^3$\,${\rm S/m}^{-2}$, which corresponds to good contacts in the framework of the smoothened CEM; see~\cite[Section~4.3]{Hyvonen17b}. For the lagged diffusivity iteration with the smoothened TV prior in Section~\ref{sec:TV} this means that the linearization of the forward model is always performed at the correct value of $\zeta$; in the terminology of Remark~\ref{rmrk:seq_lin}, the update formula \eqref{eq:xi_update} is not utilized, but one simply sets $\xi_k = 0$. Although it is not realistic to assume the contact conductances are known in a practical measurement setup for EIT, we have decided to reduce their role in our numerical tests to enable concentrating on the effect the head geometry and the location of the ROI or a stroke have on the A-optimized electrode positions. The conductivities of healthy skin, skull and brain tissues are chosen as $0.2$\,S/m, $0.06$\,S/m and $0.2$\,S/m, respectively (cf.~\cite{Gabriel96,Lai05,Latikka01,McCann19,Oostendorp00}).

The zero-mean Gaussian measurement noise is assumed to have a diagonal covariance matrix of the form $\eta^2 {\rm I}$, where the standard deviation is defined to be
\begin{equation}
  \label{eq:noise_std}
\eta = \omega \big( \max \, \mathcal{U}(\sigma,\zeta; I_1, \dots, I_{M-1}) -  \min \, \mathcal{U}(\sigma,\zeta; I_1, \dots, I_{M-1}) \big).
\end{equation}
Here, max and min operations are taken over the components of the complete measurement vector $\mathcal{U}\in \R^{M(M-1)}$, the conductivity $\sigma$ is defined by the background levels for the three-layer head model, and the current patterns are those for the symmetric initial electrode configuration shown on the first row of Figure~\ref{fig:init_conf}. We choose the value $\omega = 10^{-3}$ for the free scaling parameter and keep the noise covariance matrix fixed throughout the numerical studies.

The free parameters in Algorithms~\ref{alg:GD} and \ref{alg:armijo_alg} are set to $\epsilon = 0$, $N_{\rm GD} = 40$, $\lambda = 0.5$, $N_{\mathrm{Armijo}} = 30$, $\alpha = 1/2$ and $\beta = 5/6$. In particular, we always simply run Algorithm~\ref{alg:GD} for 40 rounds without considering stopping it earlier because the effect of numerical errors makes it difficult to reliably infer whether a local minimum has been reached. The choices of $\epsilon$, $\lambda$, $\alpha$ and $\beta$ have an effect on the performance of the optimization routine, but this aspect is not further analyzed here.

\subsection{Gaussian prior with a region of interest}
\label{sec:gauss_numer}
Let us start with a Gaussian prior for the conductivity perturbation, which allows offline minimization for the A-optimality target $\Psi_{\rm A}$ in our linearized setting. The prior covariance is given elementwise as 
\begin{equation}
 \label{eq:prior_cov}
 (\Gamma_{\rm prior})_{i,j} = \varsigma^2 \exp \left(-\frac{| x_i - x_j |^2}{2\ell^2} \right),
 \end{equation}
where $\ell = 0.05$ is the assumed correlation length, $\varsigma = 0.2$ is the pixelwise standard deviation, and $x_i$ and $x_j$ are the coordinates of the nodes with indices $i$ and $j$ in the FE mesh.

\begin{figure}[t]
\center{
      {\includegraphics[width=6cm]{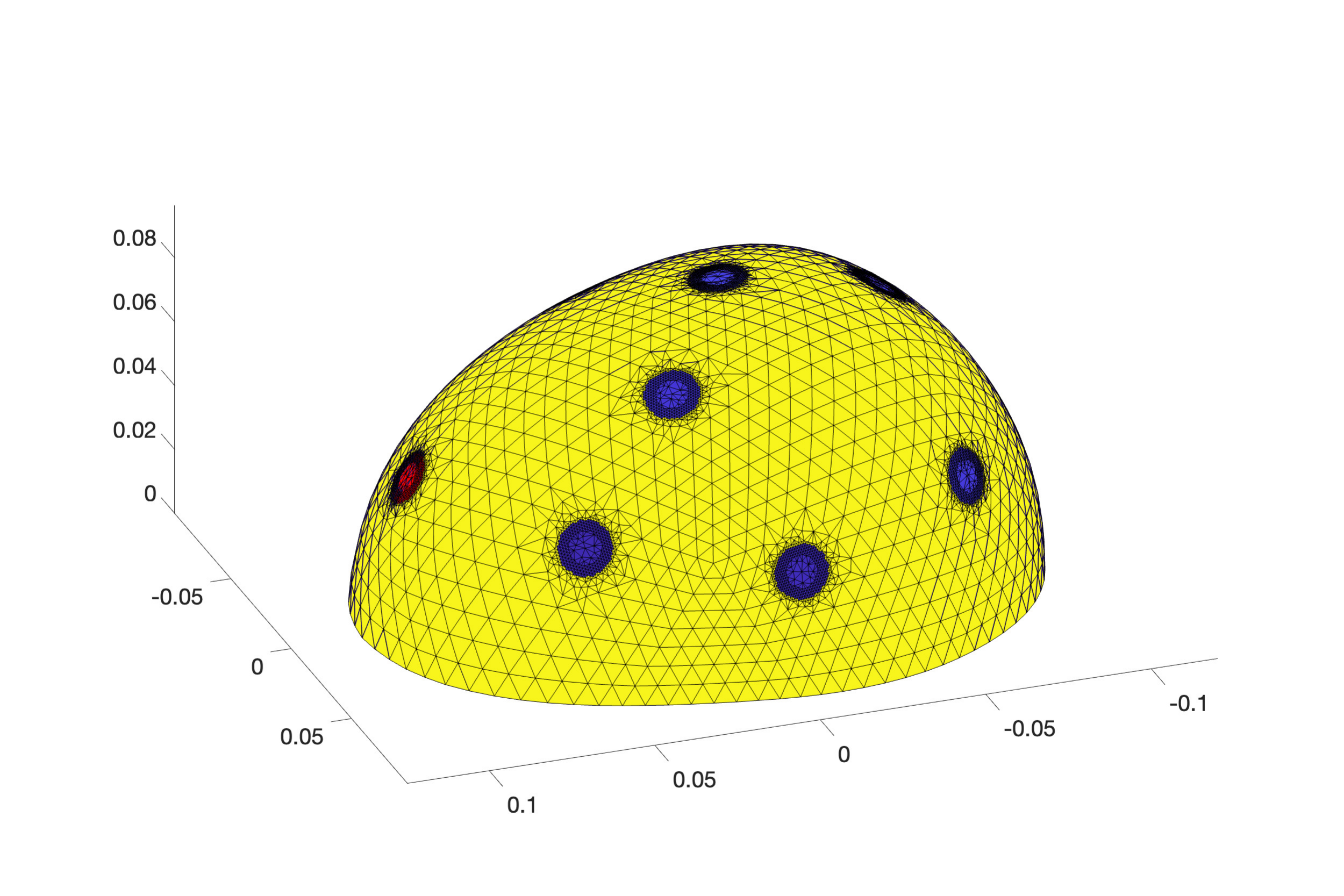}}
  \quad     {\includegraphics[width=6cm]{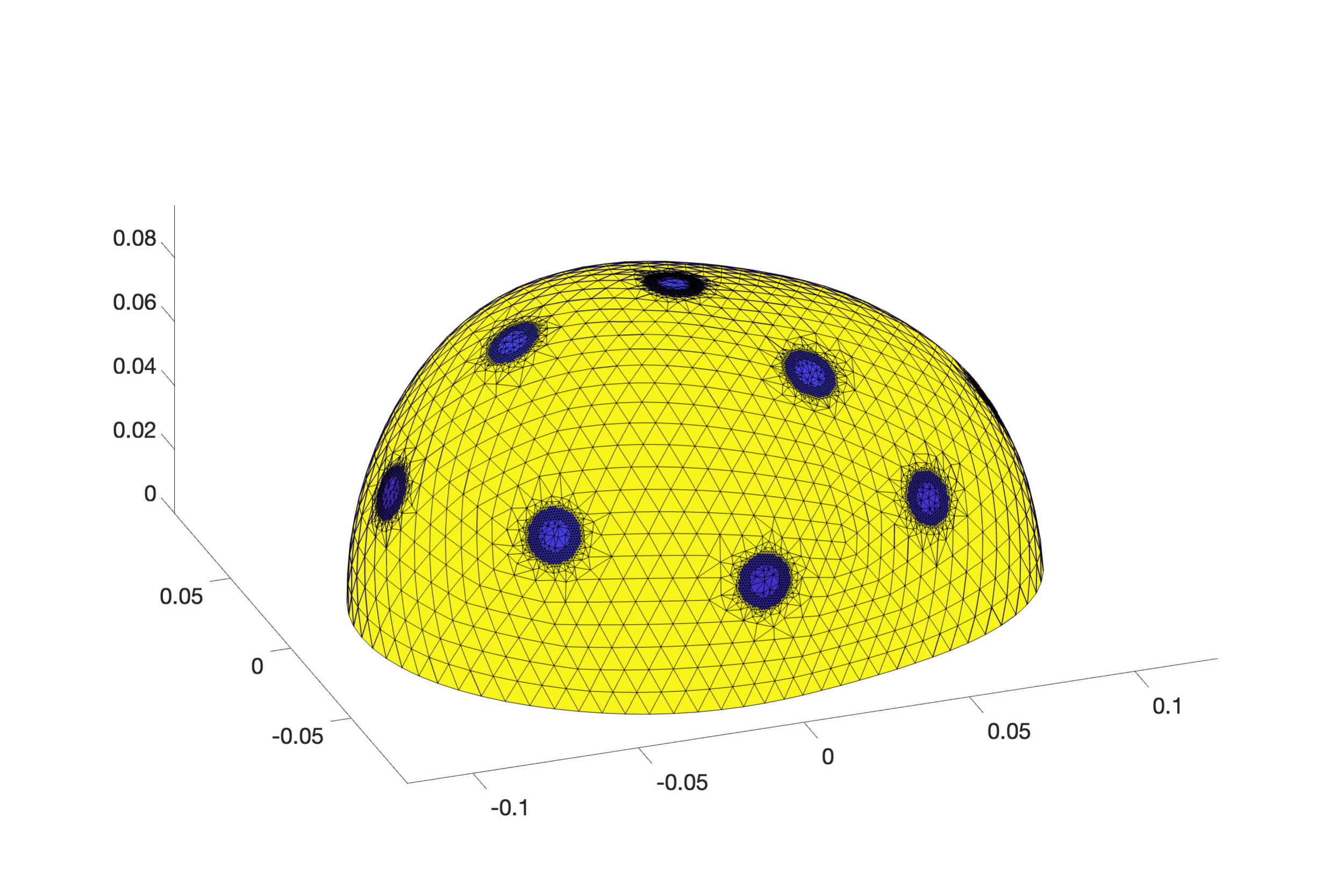}}
      \linebreak
       {\includegraphics[width=5.5cm]{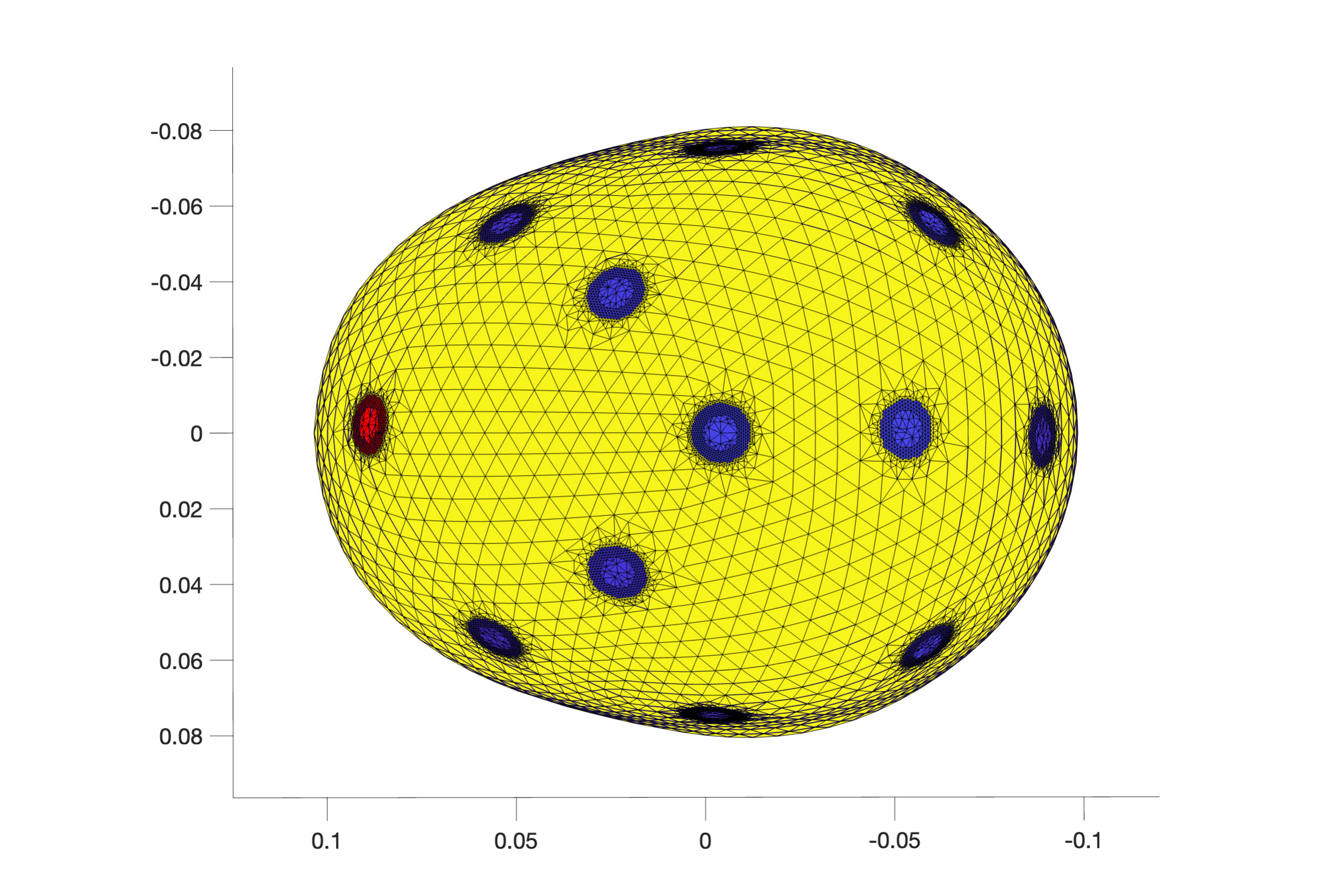}}  \quad
      {\includegraphics[width=5.5cm]{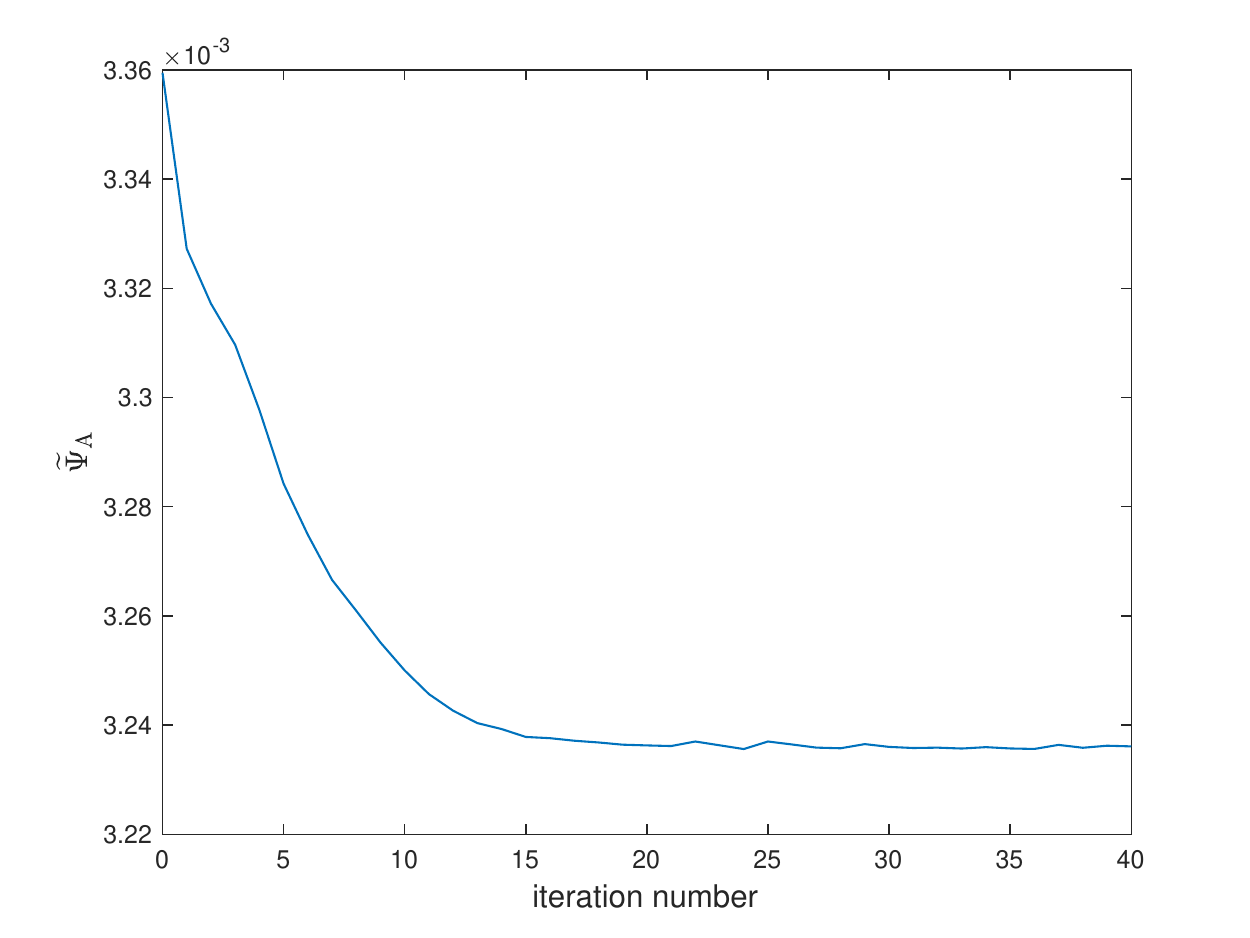}}
  }
\caption{Gaussian Test 1. The ROI is the brain with its bottom layer of height $2$\,{\rm cm} excluded. Top and bottom left: the optimized electrode setup deduced by starting from the symmetric initial configuration of Figure~\ref{fig:init_conf}. Bottom right: the evolution of $\widetilde{\Psi}_{\rm A}$ during the application of Algorithm~\ref{alg:GD}.
}
\label{fig:Gauss_test1}
\end{figure}

The first numerical test assumes the ROI is the intersection of the brain with the half-space $\{ x \in \R^3 \, | \, x_3 \geq 0.02 \}$. The initial and optimized electrode configuration are shown on the top row of Figure~\ref{fig:init_conf} and in Figure~\ref{fig:Gauss_test1}, respectively. The electrodes have reorganized themselves but not moved much in total, which is not surprising since the initial measurement configuration and the ROI are both symmetric. Be that as it may, the evolution of $\widetilde{\Psi}_{\rm A}$ plotted in the bottom right image of~Figure~\ref{fig:Gauss_test1} reveals a decrease of $3.67$\% during the minimization process. If the feeding electrode in the initial configuration on the top row of Figure~\ref{fig:init_conf} and in the setup of Figure~\ref{fig:Gauss_test1} is chosen to be another electrode, the reduction in the value of $\widetilde{\Psi}_{\rm A}$ when switching from the initial configuration to the one in Figure~\ref{fig:Gauss_test1} is still $3.14$--$3.80$\% depending on the chosen feeding electrode. Take note that the optimization algorithm is not rerun to deduce these numbers, but the index of the current-feeding electrode is simply chosen to go through all other options when evaluating  $\widetilde{\Psi}_{\rm A}$ for the two electrode configurations. Although the choice of the feeding electrode has an effect on the A-optimal configuration, optimizing the configuration for a fixed feeding electrode thus seems to reduce the value of the A-optimality target for other choices of the feeding electrode as well --- at least for the considered parameter values.

\begin{figure}[t]
\center{
  {\includegraphics[width=6cm]{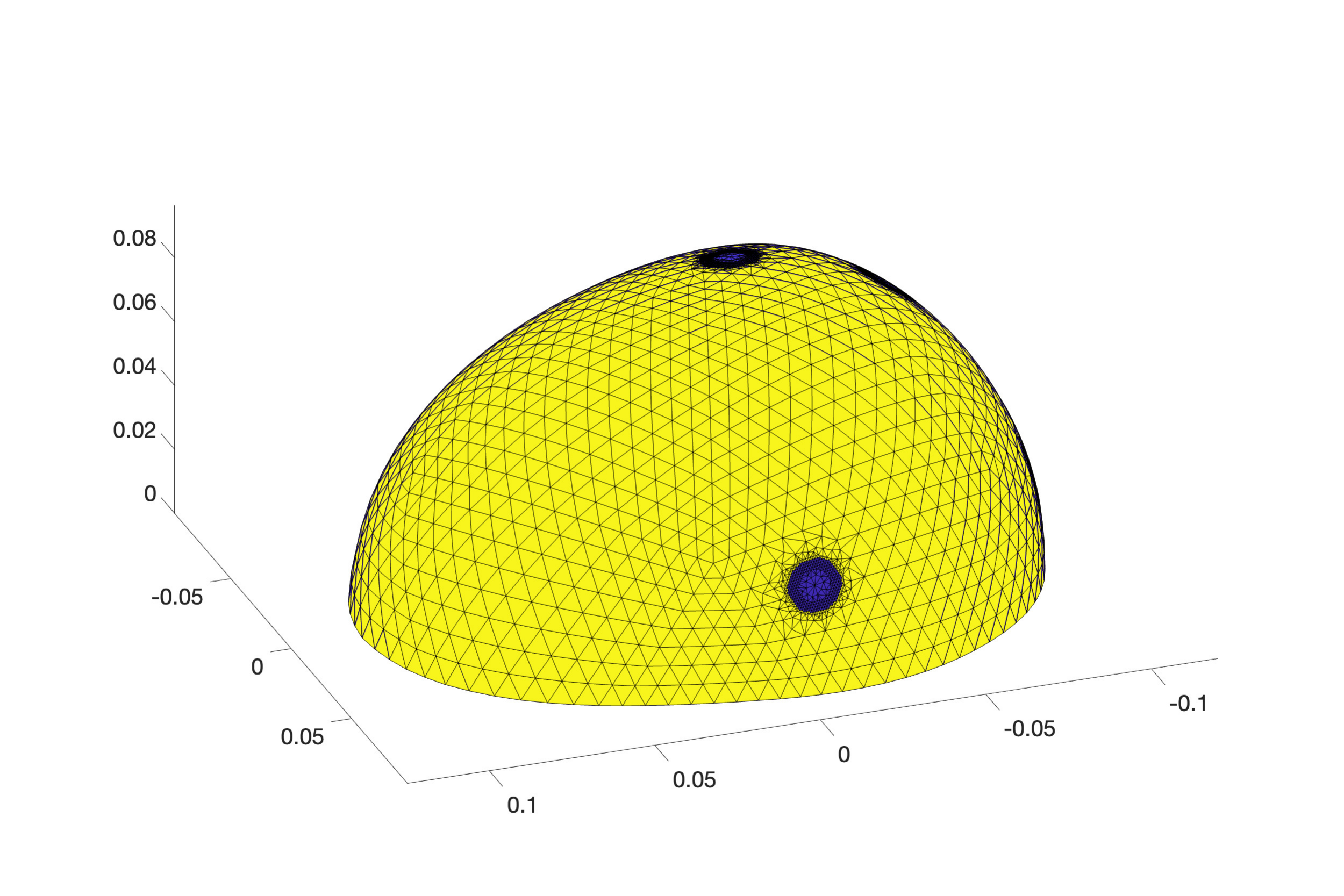}}
  \quad
      {\includegraphics[width=6cm]{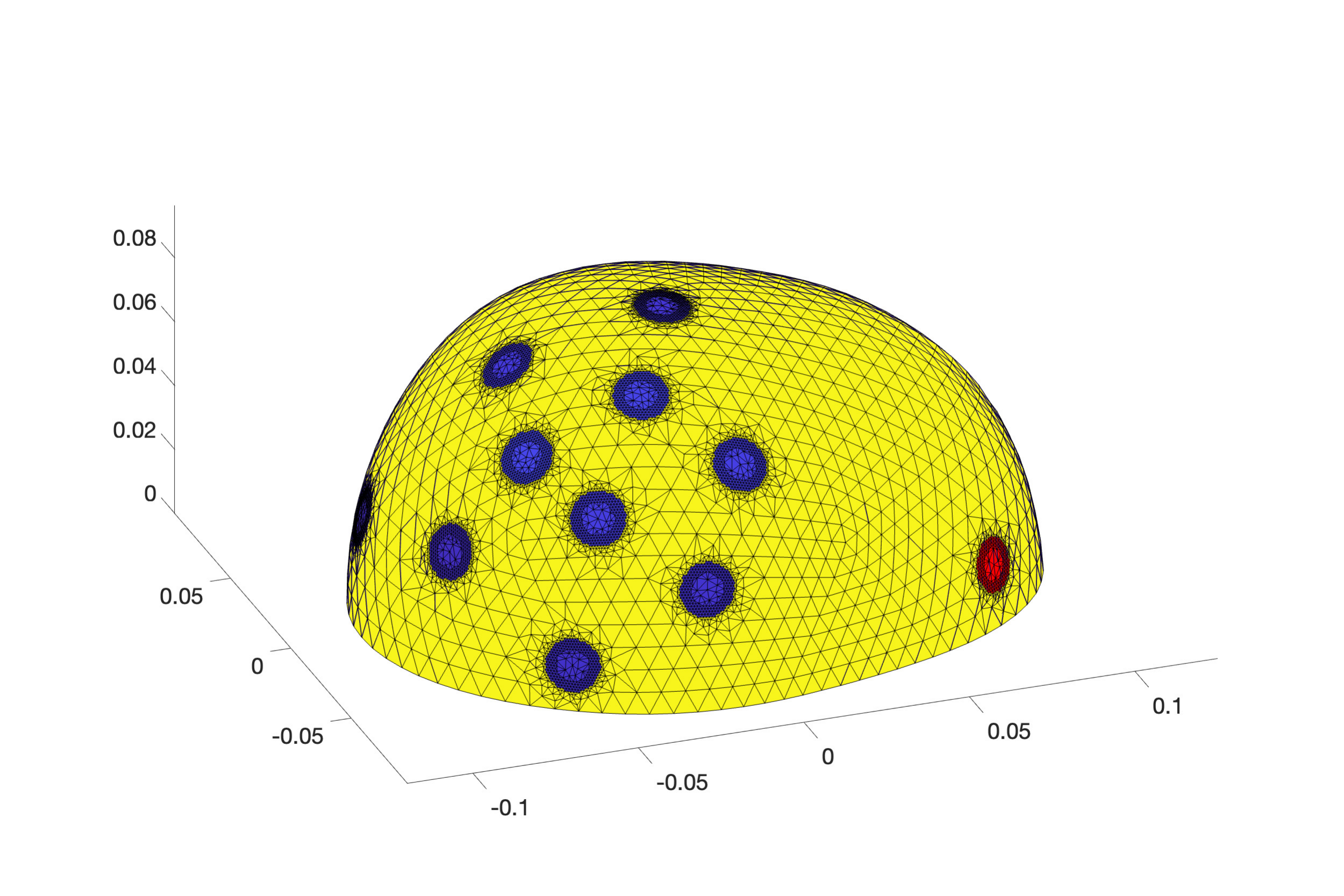}} \linebreak
      {\includegraphics[width=6cm]{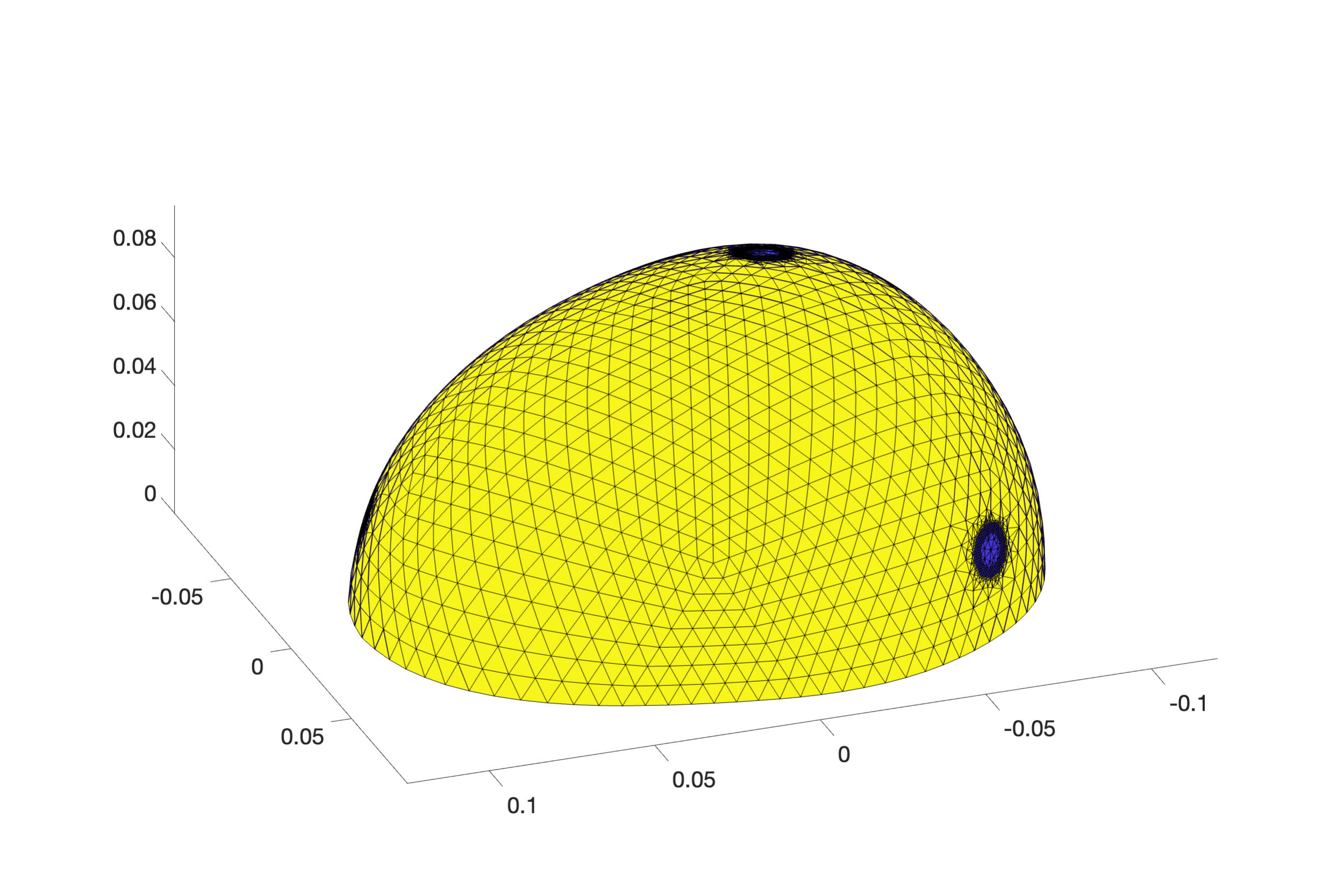}} \quad {\includegraphics[width=6cm]{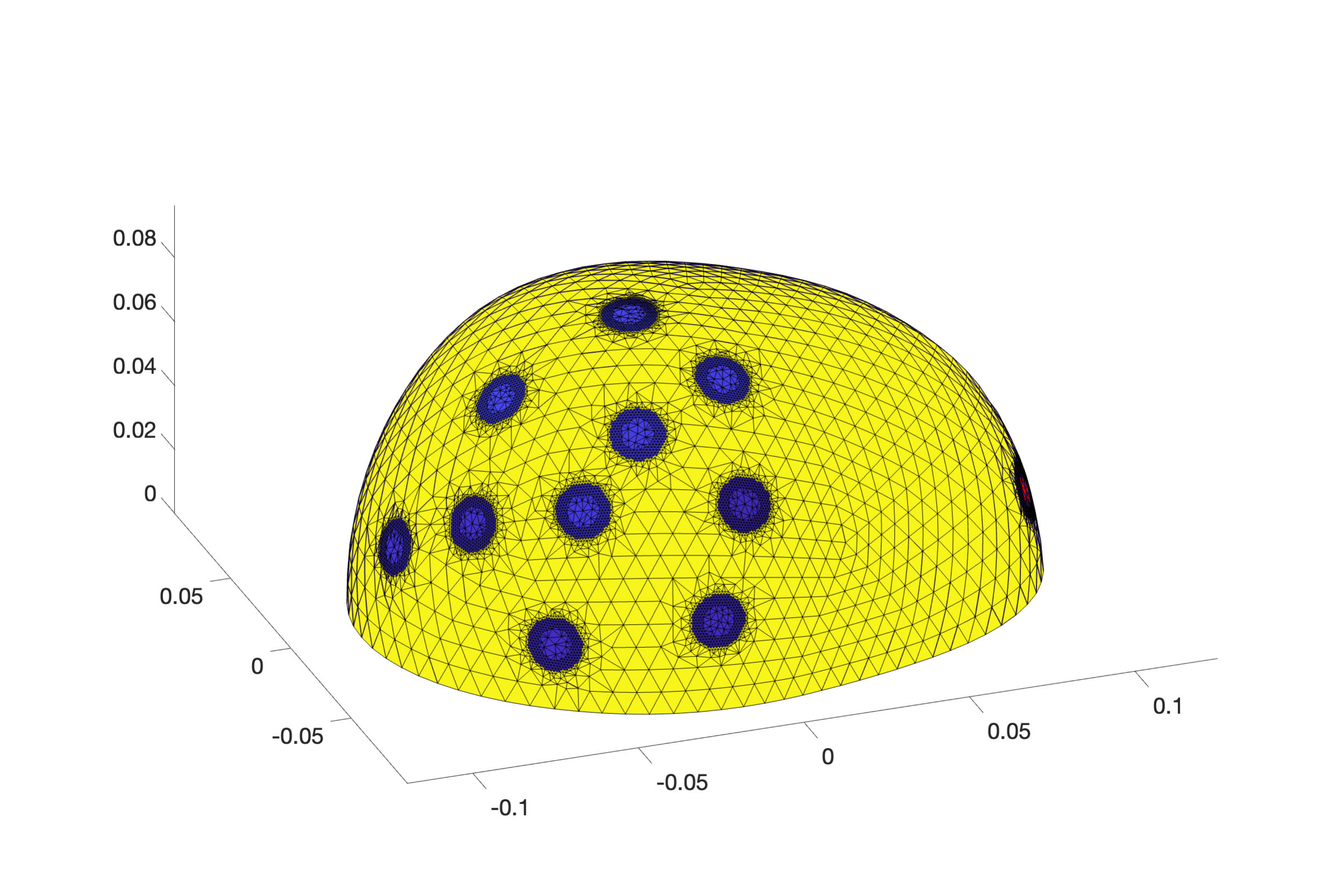}}
  }
\caption{Gaussian Test 2. The ROI is the right posterior quadrant of the brain with its bottom layer of height $2$\,{\rm cm} excluded. Top: the optimized electrode setup deduced by starting from the symmetric initial configuration on the top row of Figure~\ref{fig:init_conf}. Bottom: the optimized electrode setup deduced by starting from the pre-optimized initial configuration on the bottom row of Figure~\ref{fig:init_conf}.
}
\label{fig:Gauss_test2_1}
\end{figure}

In the second test the ROI is chosen to be the intersection of the brain with the three half spaces
$$
\{ x \in \R^3 \, | \, x_3 \geq 0.02  \}, \qquad \{ x \in \R^3 \, | \, x_1 \leq 0 \}, \qquad \{ x \in \R^3 \, | \, x_2 \leq 0 \}.
$$
In other words, it is essentially the right posterior quadrant of the brain, which could,~e.g.,~reflect one's knowledge on the position of the stroke based on a CT image. The top and bottom rows of Figure~\ref{fig:Gauss_test2_1} illustrate the optimized electrode positions starting from the top and bottom initial configurations in Figure~\ref{fig:init_conf}, respectively. For both initial configurations, the minimization algorithm moves most of the electrodes to the vicinity of the right posterior quadrant of the brain, with a couple of electrodes, including the feeding one, situated at a longer distance from the ROI. Although qualitatively similar, the optimized configurations are not precisely the same, which may be caused by local minima or numerical inaccuracies in the evaluation of $\Psi_{\rm A}$ and its derivatives during the minimization process.

\begin{figure}[t]
\center{
  {\includegraphics[width=5.5cm]{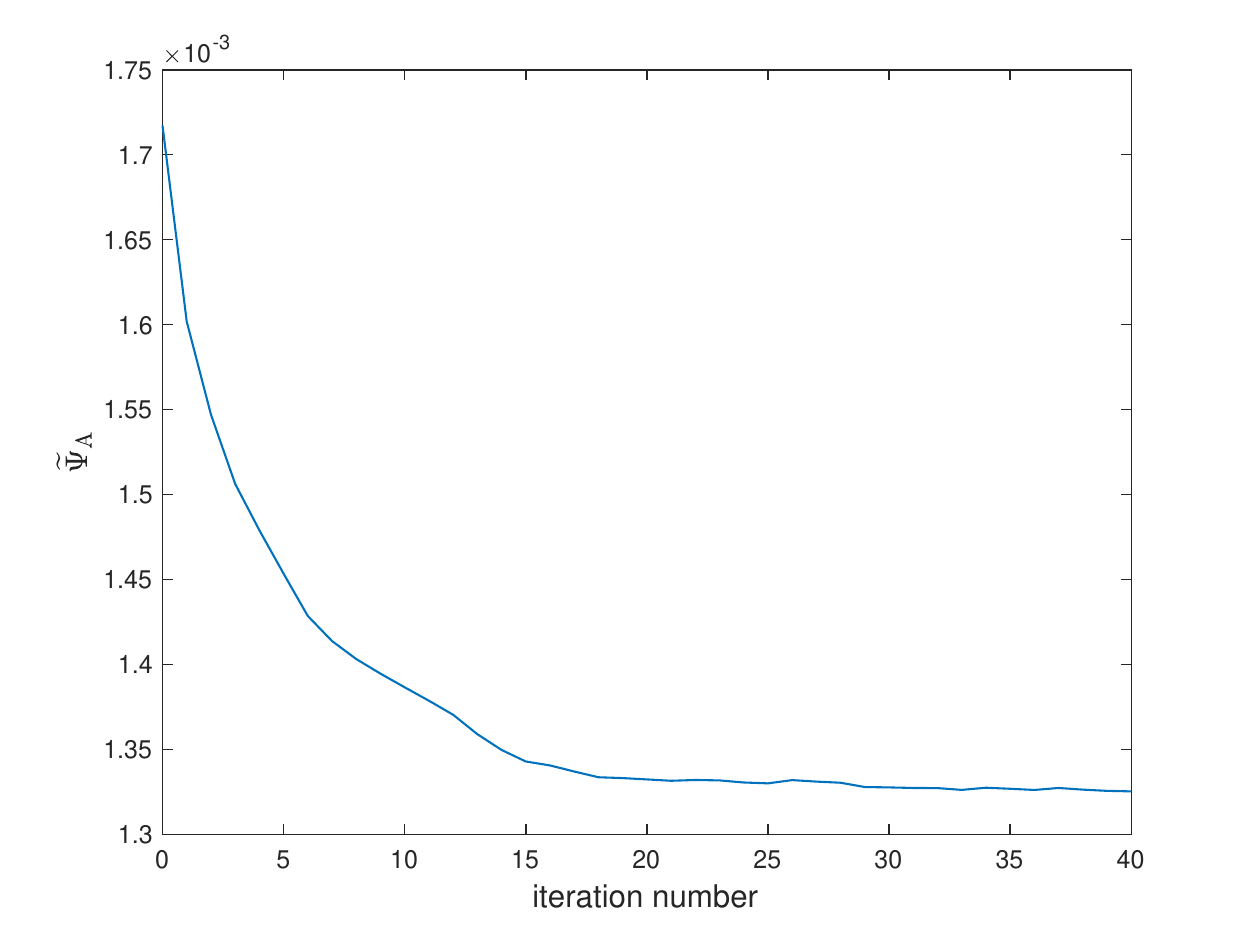}}
  \quad
      {\includegraphics[width=5.5cm]{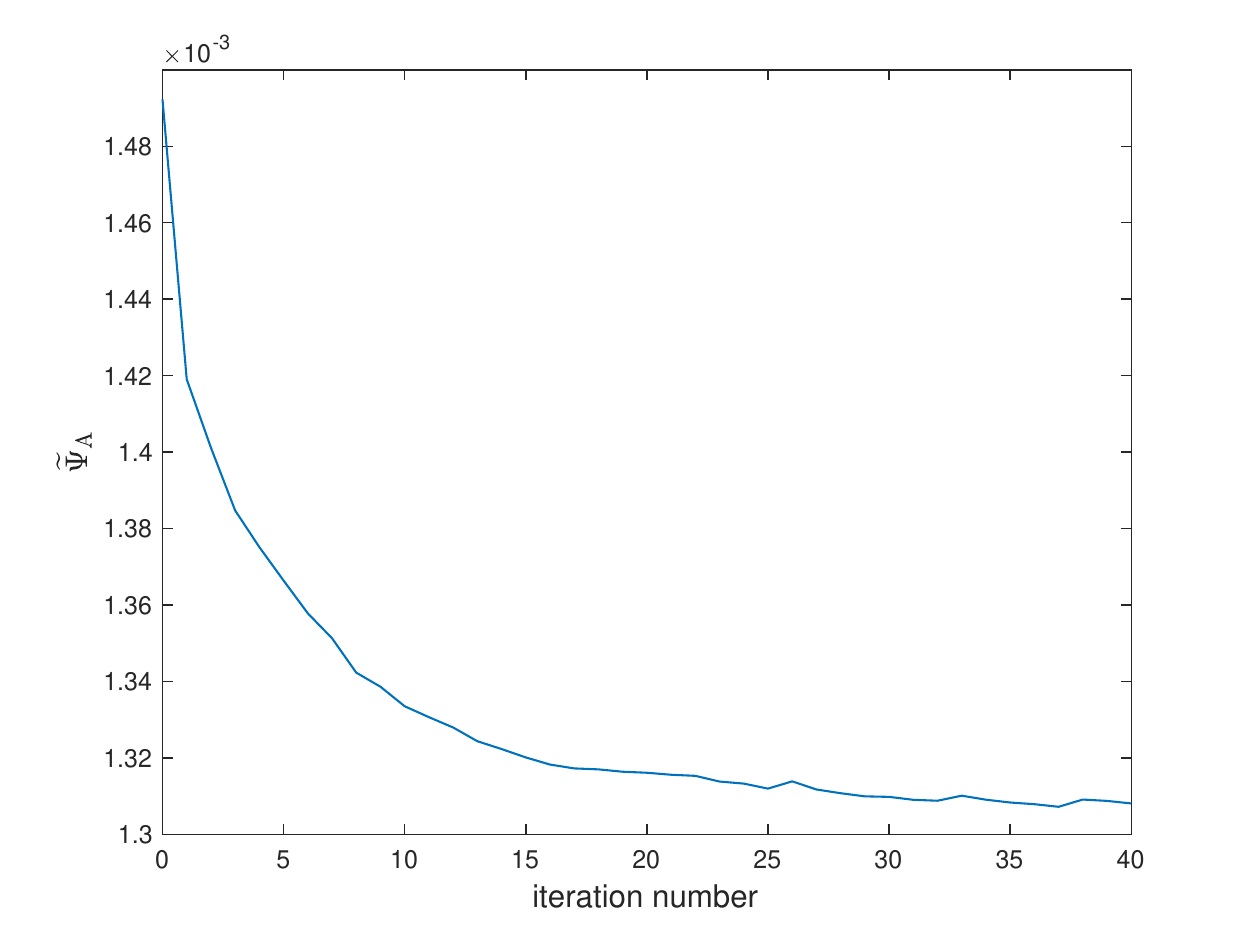}}
  }
\caption{Gaussian Test 2. The evolution of $\widetilde{\Psi}_{\rm A}$ during Algorithm~\ref{alg:GD} when the ROI is the right posterior quadrant of the brain. Left: the symmetric initial configuration on the top row of Figure~\ref{fig:init_conf}. Right: the pre-optimized initial configuration on the bottom row of Figure~\ref{fig:init_conf}.
}
\label{fig:Gauss_test2_2}
\end{figure}

Figure~\ref{fig:Gauss_test2_2} shows the evolution of $\widetilde{\Psi}_{\rm A}$ during Algorithm~\ref{alg:GD} for the considered two initial configurations. The initial electrode setup that is concentrated in the vicinity of the ROI leads to a slightly,~i.e.,~$1.30$\%, lower minimized value for $\widetilde{\Psi}_{\rm A}$. However, the difference in the corresponding initial values for $\widetilde{\Psi}_{\rm A}$ is significant: the square root of the expected squared $L^2$ reconstruction error over the ROI for the symmetric initial configuration is $15.08$\% higher than that for the initial configuration accounting for the location of the ROI to a certain extent. This demonstrates that the intuition of placing most available electrodes close to the ROI seems correct. For the symmetric initial configuration, running Algorithm~\ref{alg:GD} resulted in a reduction of $22.82$\% in the value of~$\widetilde{\Psi}_{\rm A}$, whereas the reduction was only $12.34$\% for the `pre-optimized' initial configuration.

\subsection{Total variation with an inclusion}
\label{sec:TV_numer}

\begin{figure}[t]
\center{
  {\includegraphics[width=12cm]{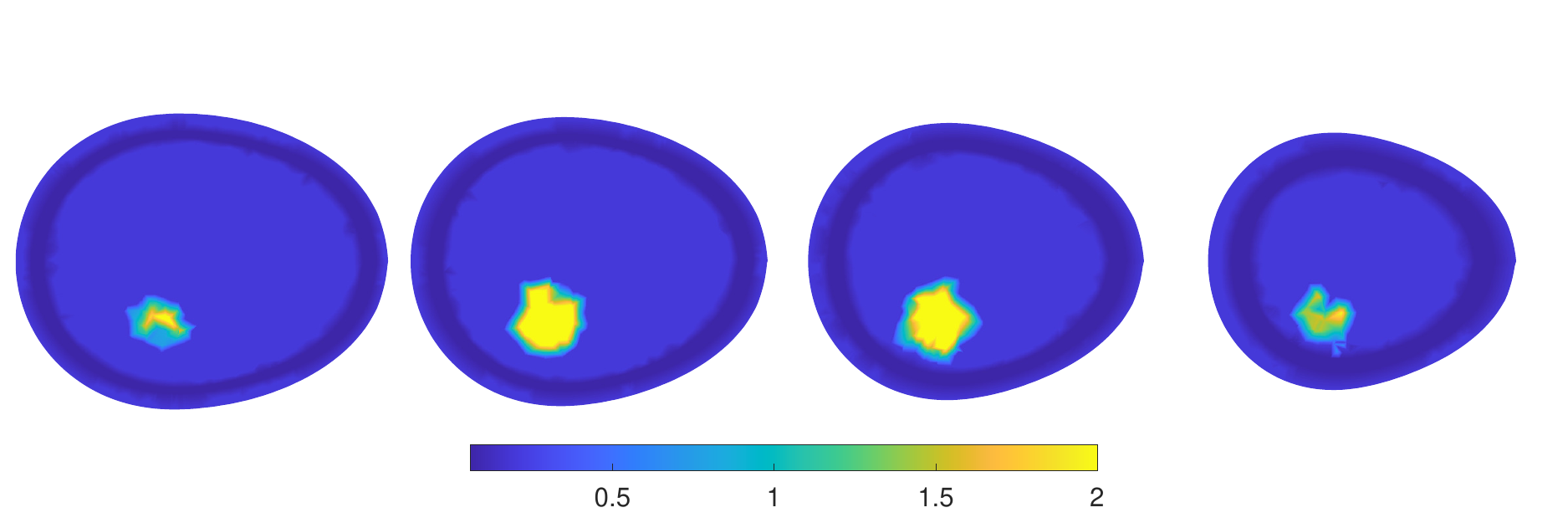}}
}
\caption{Cross-sections of the investigated head model with a spherical conductive inclusion in the right posterior quadrant of the brain. The horizontal slices are at heights $3$\,cm, $4$\,cm, $5$\,cm and $6$\,cm. The unit of conductivity is S/m.
  }
\label{fig:inclusion}
\end{figure}

Assume there is a spherical conductive anomaly of radius $1.8$\,cm in the right posterior quadrant of the brain as illustrated in Figure~\ref{fig:inclusion}. The inhomogeneity could model a hemorrhagic stroke \cite{Toivanen21}. We first measure (i.e.,~simulate) noisy data (cf.~\eqref{eq:noise_std}) with an initial electrode configuration and then compute a reconstruction corresponding to a smoothened TV prior by combining lagged diffusivity iteration with sequential linearizations, as explained in Section~\ref{sec:TV}. The lagged diffusivity iteration associated to the final linearization of the measurement model can be interpreted to produce a Gaussian posterior for the conductivity perturbation, with its mean and covariance given by \eqref{eq:LD_mean} and \eqref{eq:LD_covariance}, respectively, cf.~Remark~\ref{rmrk:seq_lin}. This covariance is finally used as $\Gamma_{\rm prior}$ in \eqref{eq:Gaussian_posterior} and indirectly in \eqref{eq:Aoptimal} to enable finding new, more informative positions for the available electrodes via minimizing the induced A-optimality target.

In the following numerical tests, we choose the values $T = 10^{-6}$, $\gamma = 10^5$, $c_\upsilon = 300$ and $b_\upsilon = 0.01$ for the free parameters in \eqref{eq:TVexp} and compute reconstructions of the conductivity perturbation by applying five lagged diffusivity steps to each of five subsequent linearizations, which results in reconstructions of the conductivity perturbation that are adequately in line with the smoothened TV prior. We do not consider choosing the number of lagged diffusivity steps or sequential linearizations in more detail here but instead refer to~\cite{Candiani21,Harhanen15,Helin23,Helin22} for more information on this aspect in different settings.

\begin{figure}[t]
\center{
  {\includegraphics[width=12cm]{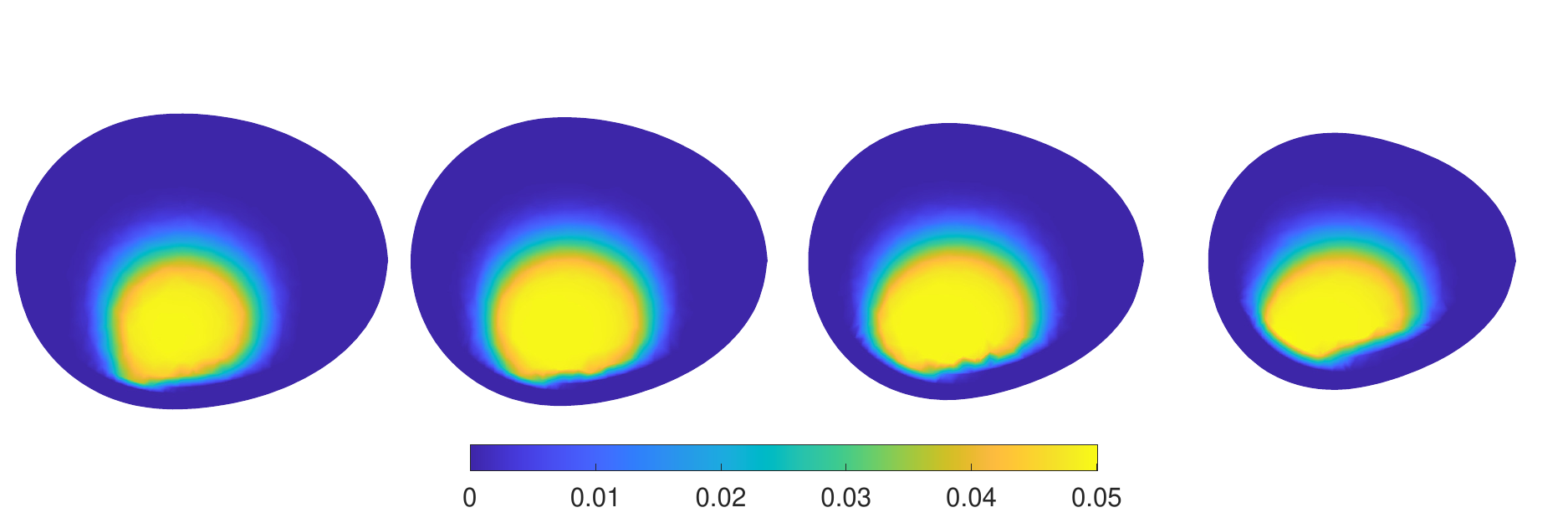}} \\
  {\includegraphics[width=12cm]{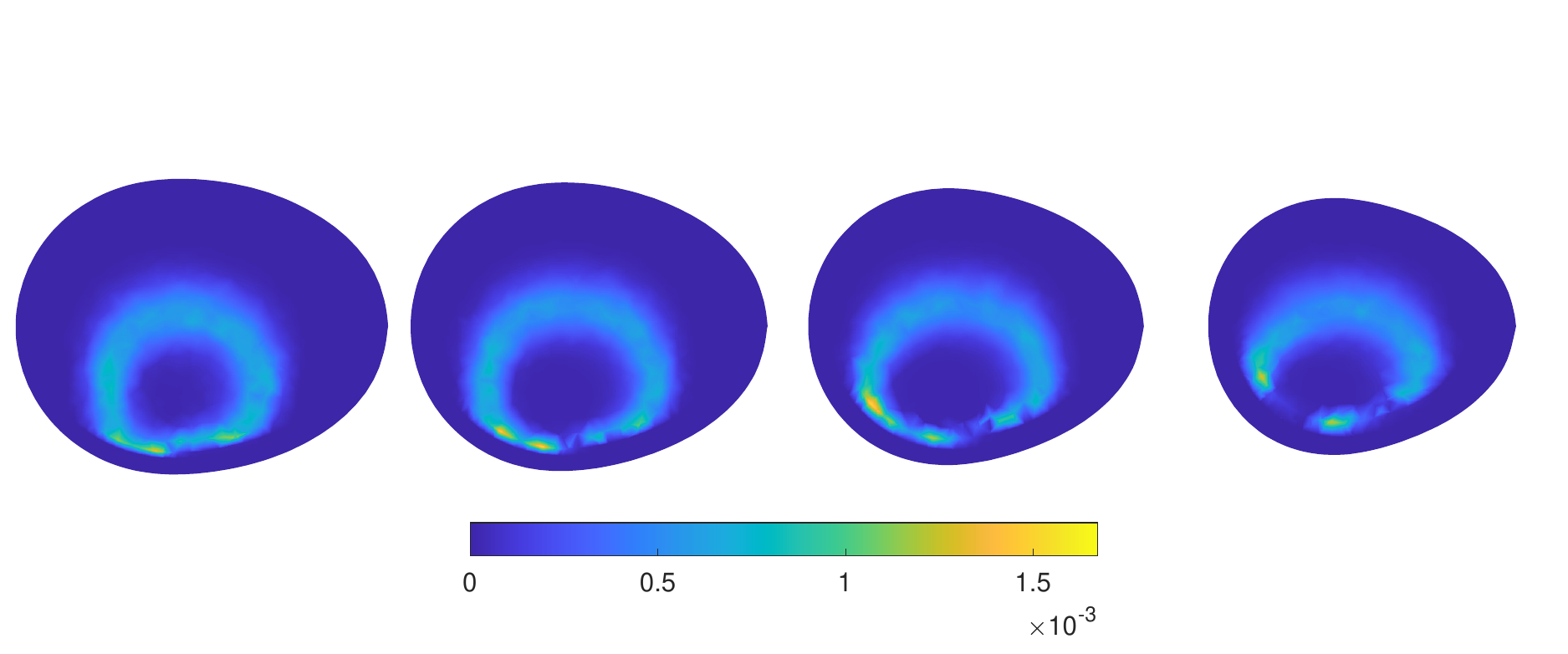}}
}
\caption{Top: cross-sections of the reconstruction of the conductivity perturbation shown in Figure~\ref{fig:inclusion} for the symmetric initial electrode configuration in Figure~\ref{fig:init_conf}. Bottom: the diagonal of the corresponding posterior covariance. The horizontal slices are at heights $3$\,cm, $4$\,cm, $5$\,cm and $6$\,cm.}
\label{fig:reco1}
\end{figure}

Let us first consider the symmetric initial electrode configuration presented on the top row of Figure~\ref{fig:init_conf}. The reconstruction computed from the associated measurements is presented at top in Figure~\ref{fig:reco1}, with the diagonal of the accompanying posterior covariance  visualized in the bottom images of the same figure. Although the diagonal of the posterior covariance, which is used as the prior for finding new positions for the electrodes, does not reveal all characteristics of the associated Gaussian distribution, it in any case indicates that the highest variances in the nodal degrees of freedom in the conductivity reconstruction are found close to the boundaries of the reconstructed inhomogeneity. It is thus to be expected that the minimization of $\Psi_{\rm A}$ leads to the electrodes moving to positions where more information about the corresponding area of high uncertainty can be obtained.

\begin{figure}[t]
\center{
  {\includegraphics[width=6cm]{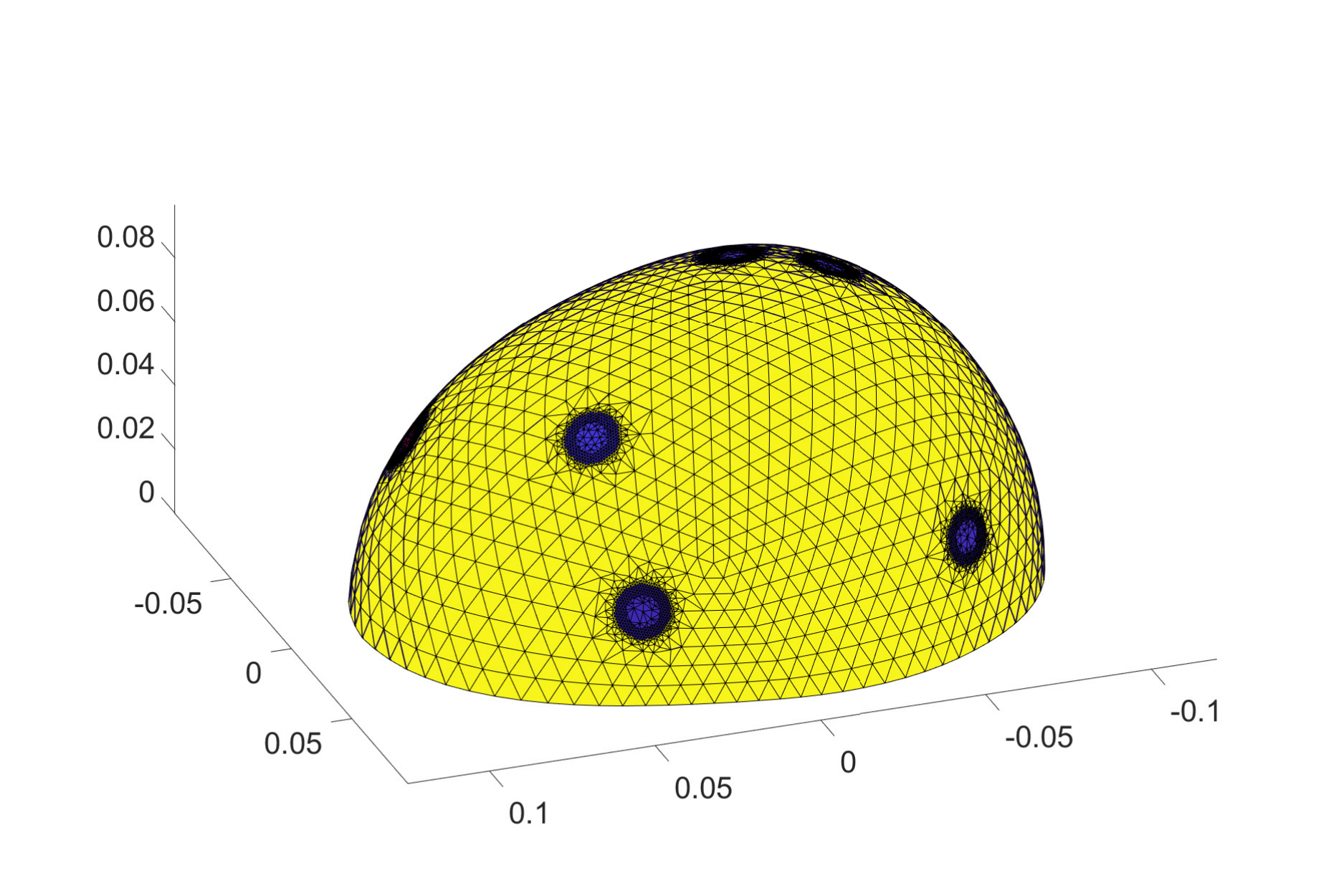}} \quad {\includegraphics[width=6cm]{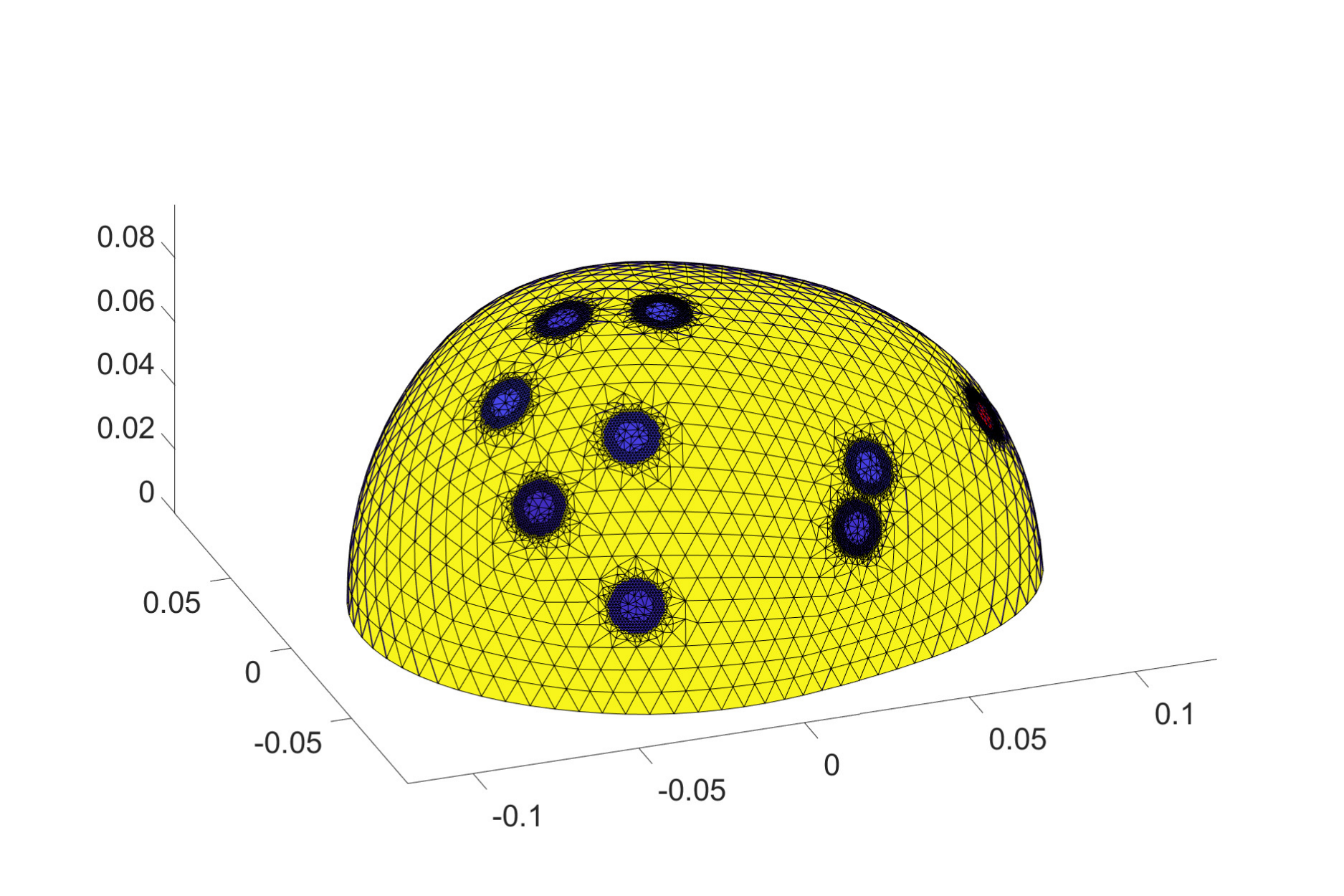}} \\
  {\includegraphics[width=6cm]{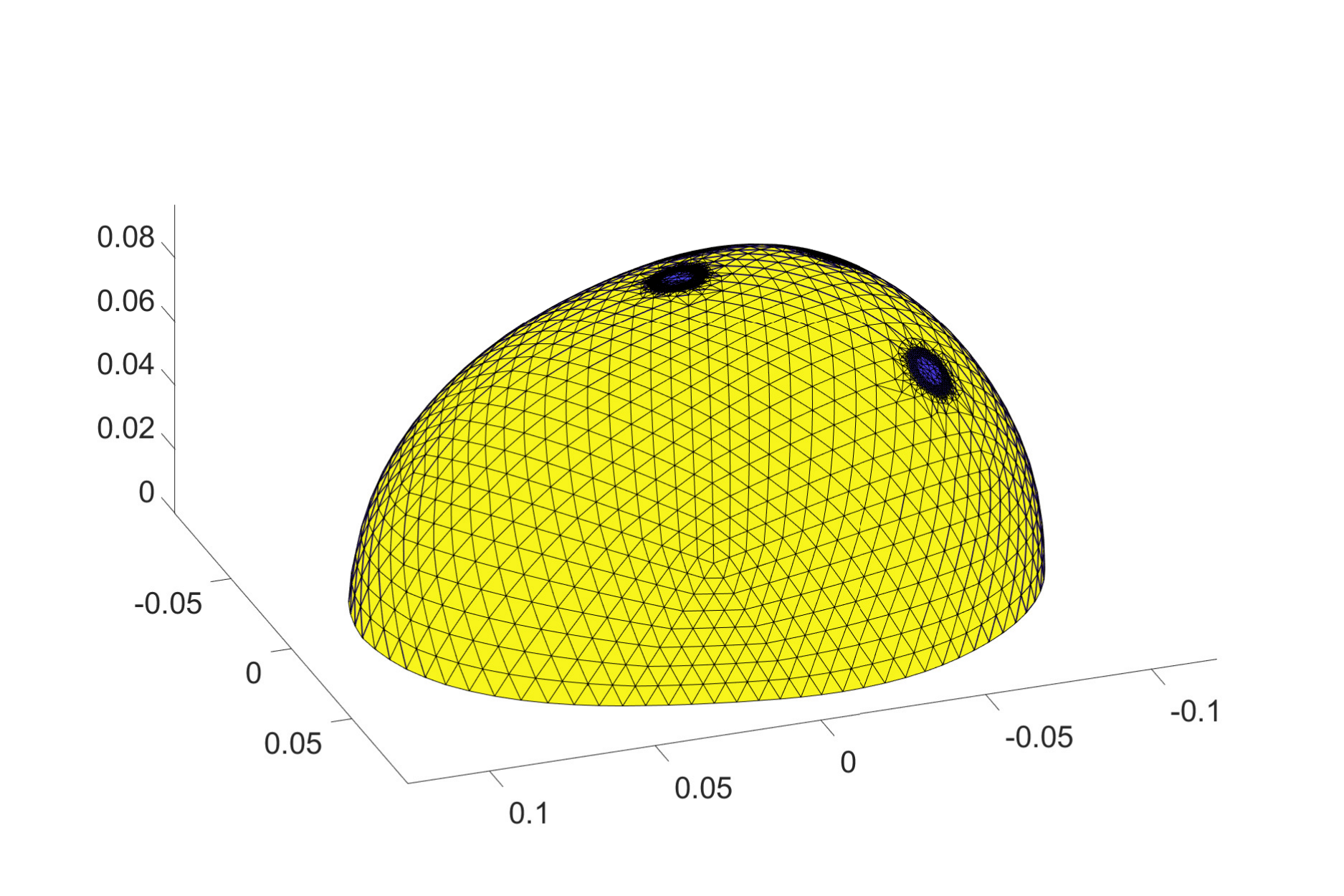}} \quad {\includegraphics[width=6cm]{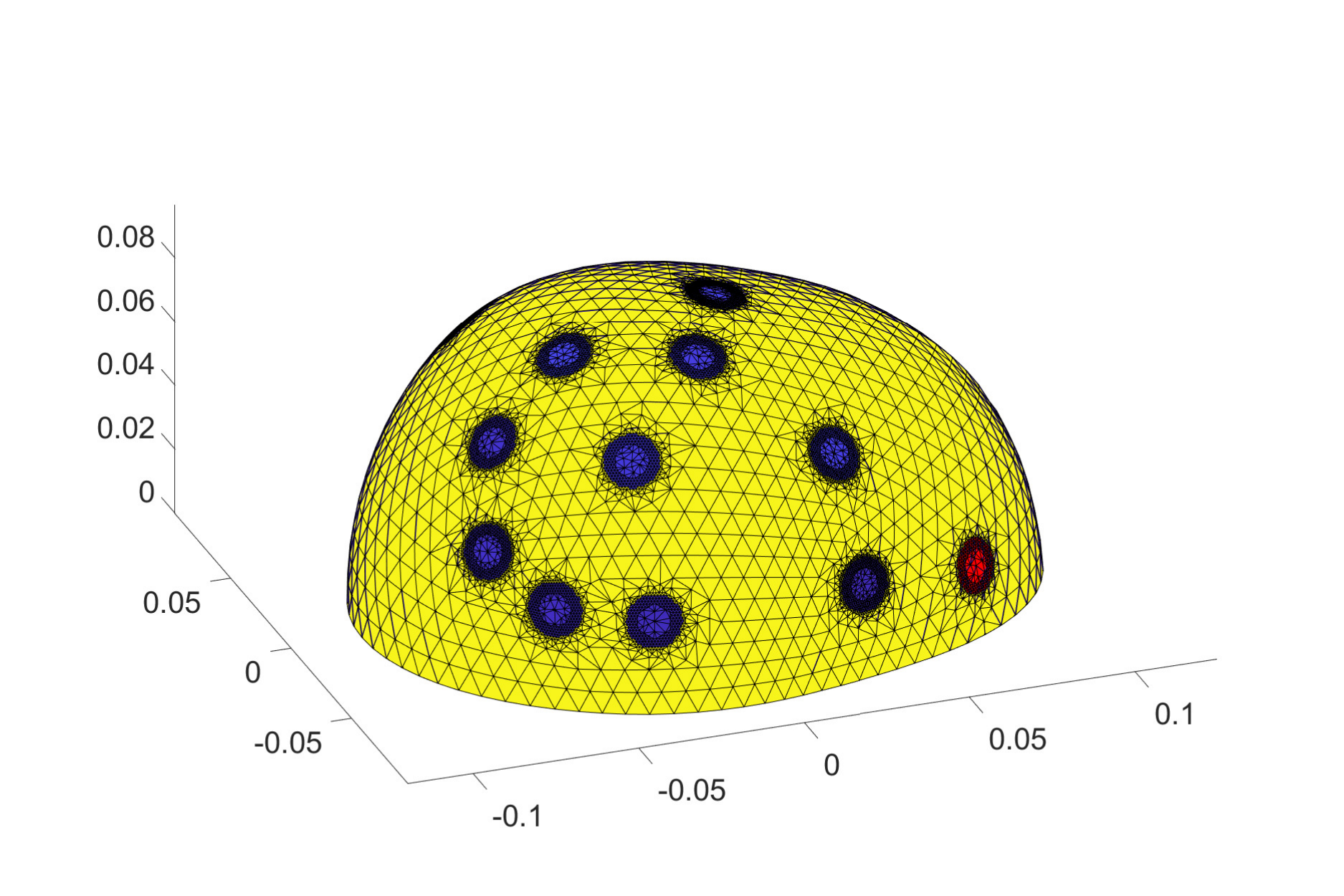}}
}
\caption{Top: the electrode configuration optimized based on the (prior) covariance obtained as a side product when forming the reconstruction in Figure~\ref{fig:reco1} from data corresponding to the symmetric electrode configuration in Figure~\ref{fig:init_conf}. Bottom: the electrode configuration optimized based on the (prior) covariance obtained as a side product when forming a reconstruction from data corresponding to the electrode configuration concentrated around the right posterior quadrant in Figure~\ref{fig:init_conf}.
  }
\label{fig:TV_test1_1}
\end{figure}

Indeed, running Algorithm~\ref{alg:GD} with the prior visualized at bottom in Figure~\ref{fig:reco1} results in the electrode positions shown in the top row of Figure~\ref{fig:TV_test1_1}. As expected, most electrodes have gathered close to the edges in the reconstruction formed from the data simulated with the initial electrode setup. Since the electrodes have heavily clustered during the optimization process, it is almost inevitable that the depicted A-optimized positions only correspond to a local minimum of $\Psi_{\rm A}$ as the electrodes presumably block one another's paths to more optimal positions.

The bottom row of Figure~\ref{fig:TV_test1_1} presents the electrode setup that results from minimizing $\Psi_{\rm A}$ with the prior covariance obtained as a side product when reconstructing the conductivity perturbation from data simulated with the electrodes concentrated around the right posterior quadrant of the head shown on the bottom row of Figure~\ref{fig:init_conf}. Since the electrodes lie initially closer to the inclusion, the obtained reconstruction of the perturbation (not shown here) is somewhat better localized to the vicinity of the true conductivity perturbation. This leads to the A-optimized electrode positions also being slightly closer to the location of the stroke than on the top row of Figure~\ref{fig:TV_test1_1}, where a couple of electrodes even remain on the uninteresting side of the head. 

\begin{figure}[t]
\center{
  {\includegraphics[width=6cm]{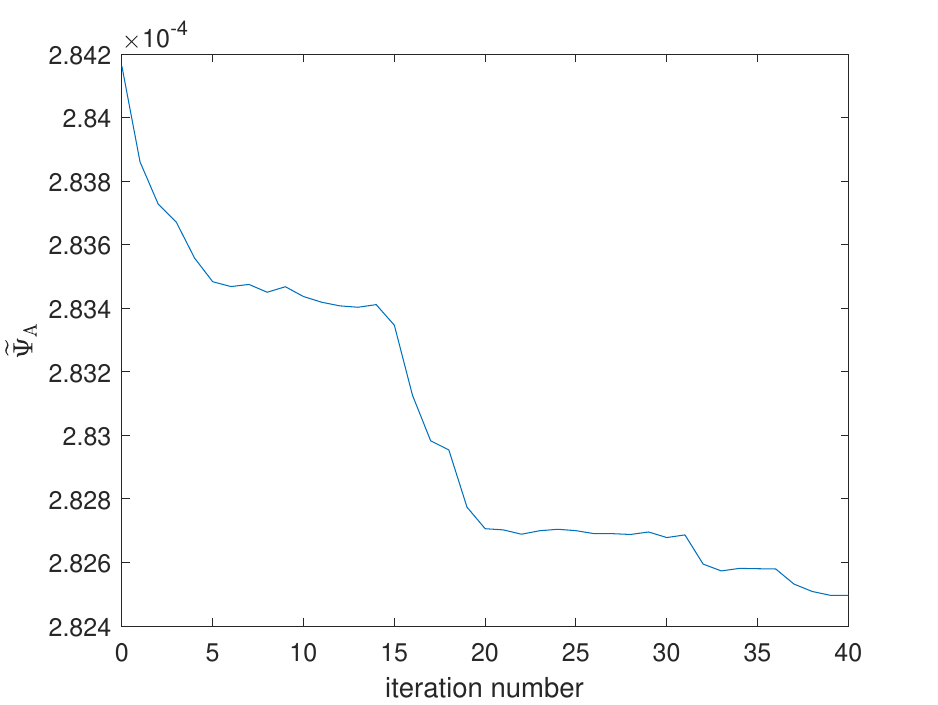}}
  \quad
      {\includegraphics[width=6cm]{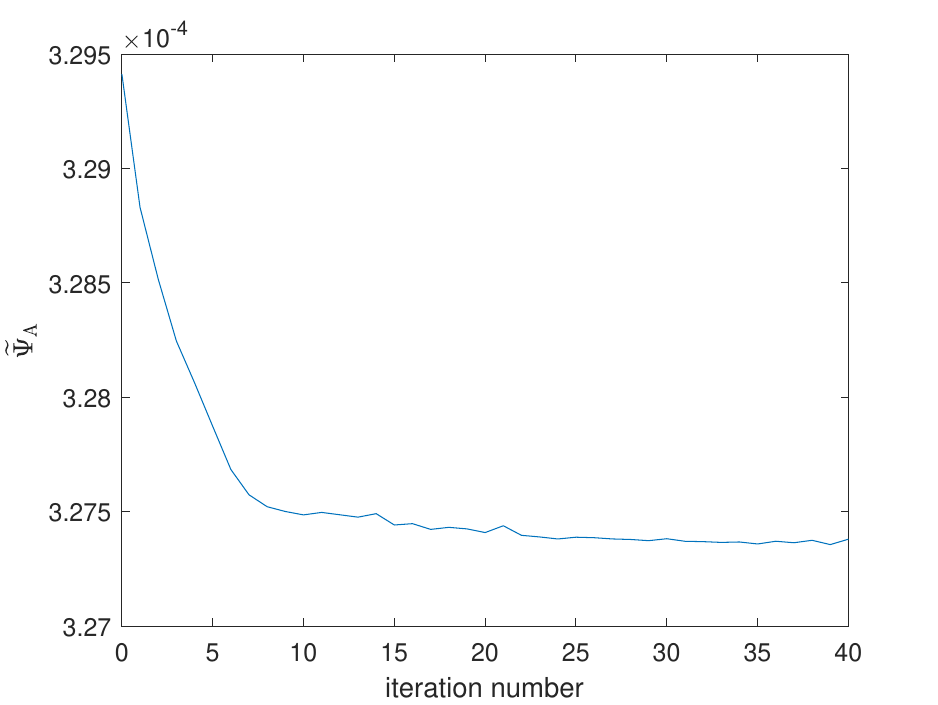}}
  }
\caption{The evolution of $\widetilde{\Psi}_{\rm A}$ during Algorithm~\ref{alg:GD}. Left: the prior covariance is obtained by forming a reconstruction from data simulated with the symmetric electrode configuration in Figure~\ref{fig:init_conf}. Right: the prior covariance is obtained by forming a reconstruction from data simulated with the  electrode configuration concentrated around the right posterior quadrant in Figure~\ref{fig:init_conf}.
}
\label{fig:TV_test1_2}
\end{figure}

The evolution of the square root of the A-optimality target, i.e.,~$\widetilde{\Psi}_{\rm A}$, for the optimization processes leading to the two electrode configurations in Figure~\ref{fig:TV_test1_1} are plotted in Figure~\ref{fig:TV_test1_2}. Note that the values of $\widetilde{\Psi}_{\rm A}$ in the two convergence plots  are not comparable as such because they correspond to two separate priors containing different levels of information on the whereabouts of the conductivity perturbation.

\begin{figure}[t]
\center{
  {\includegraphics[width=12cm]{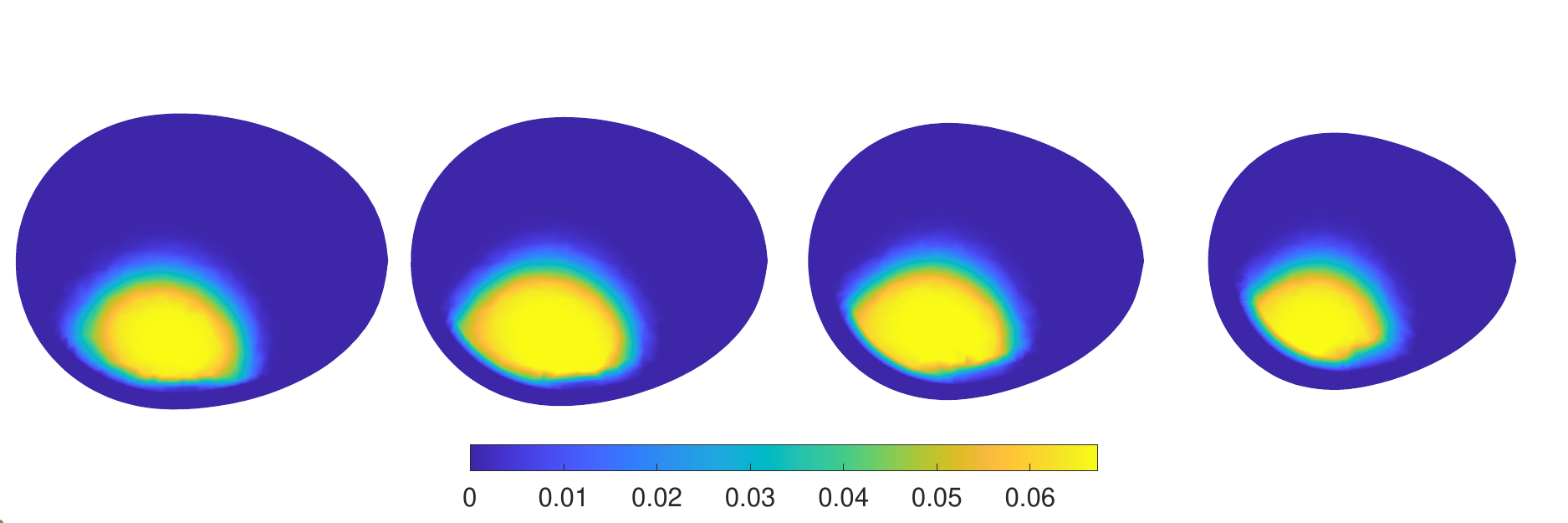}} \\
      {\includegraphics[width=6cm]{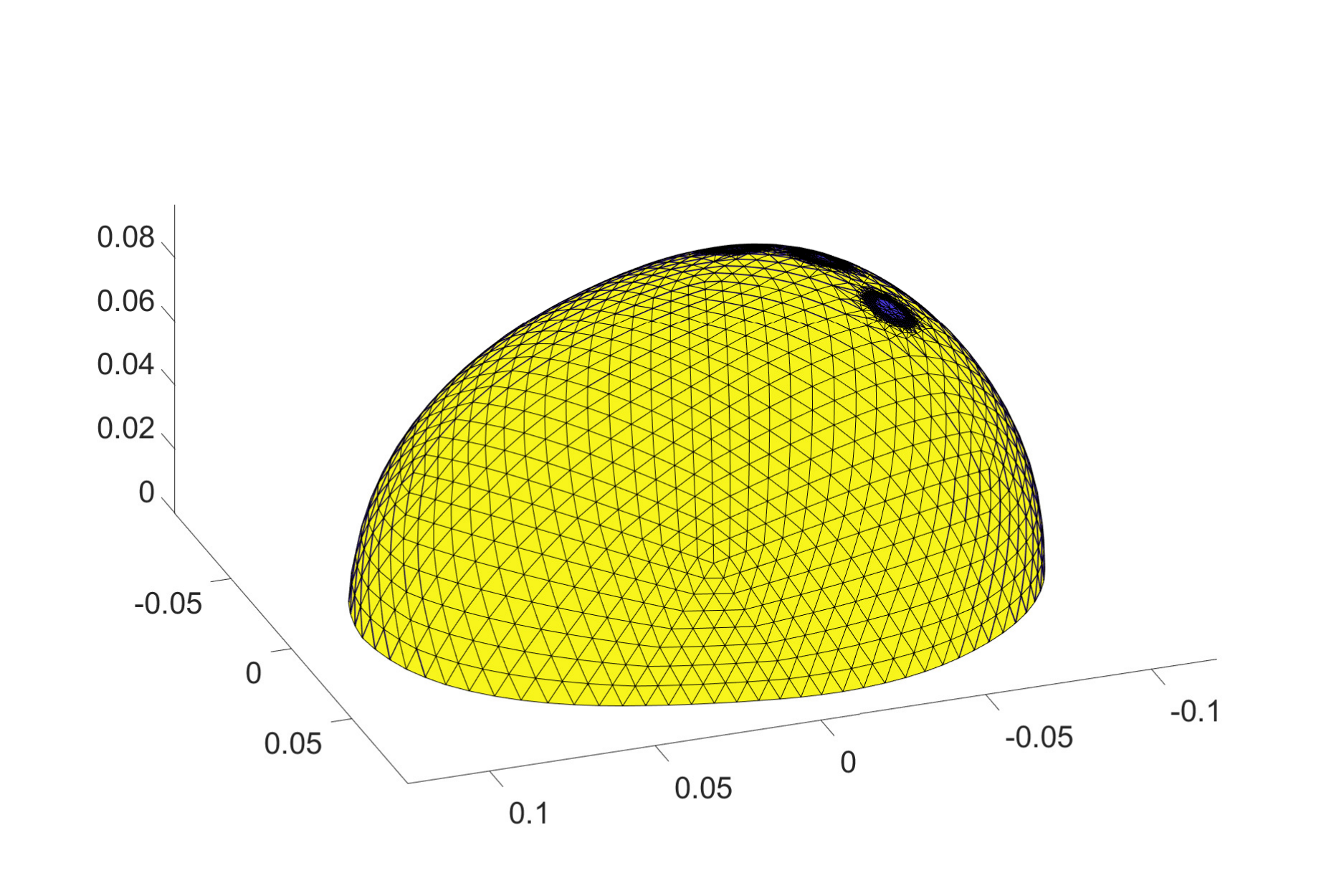}}
  \quad
      {\includegraphics[width=6cm]{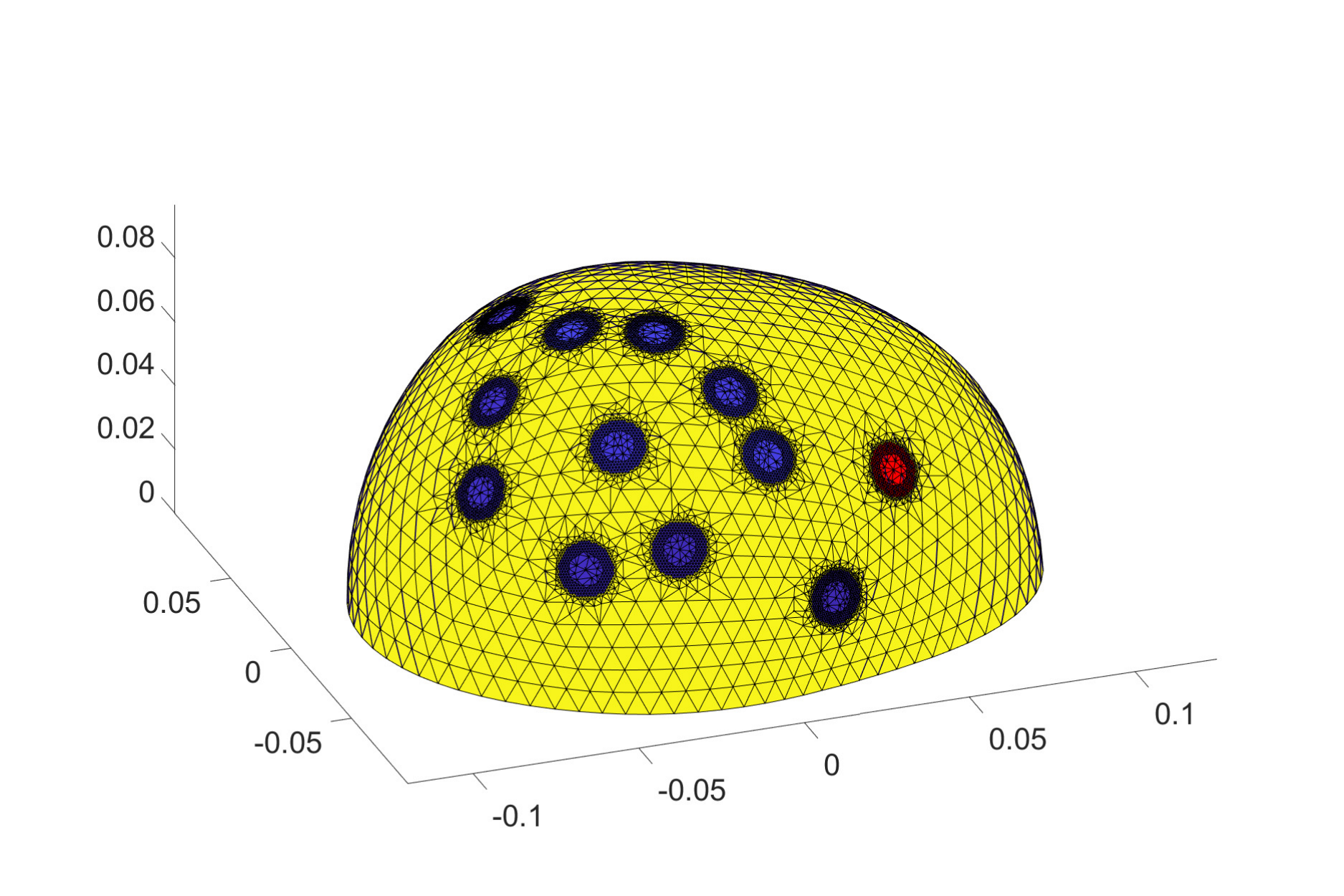}}
  }
\caption{Top: cross-sections of the reconstruction of the conductivity perturbation shown in Figure~\ref{fig:inclusion} for the electrode configuration at bottom in Figure~\ref{fig:TV_test1_1}. Bottom:  the electrode configuration optimized based on the covariance obtained as a side product when forming the presented reconstruction. 
}
\label{fig:TV_test1_3}
\end{figure}

The idea of computing a reconstruction and finding new positions for the available electrodes based on the covariance matrix formed as a side product of the reconstruction algorithm of Section~\ref{sec:TV} can be iterated. Indeed, a reconstruction computed from data simulated with the bottom electrode setup in Figure~\ref{fig:TV_test1_1} is illustrated at top in Figure~\ref{fig:TV_test1_3}, and the resulting newly A-optimized electrode setup is presented at bottom in that same figure. Since the electrodes shown in the bottom row of Figure~\ref{fig:TV_test1_1} are positioned closer to the location of the stroke than the symmetric initial electrode configuration in Figure~\ref{fig:init_conf}, the reconstruction in Figure~\ref{fig:TV_test1_3} is somewhat more accurate than that in Figure~\ref{fig:reco1}. This is also reflected in the (doubly) optimized electrodes in Figure~\ref{fig:TV_test1_3} being more tightly packed around the location of the stroke than either of the (singly) optimized configurations in Figure~\ref{fig:TV_test1_1}.

\section{Concluding remarks}
\label{sec:conclusion}

This paper tackled the optimization of electrode positions in head imaging by EIT, with the underlying motivation to monitor stroke patients in intensive care. After linearizing the smoothened CEM to introduce an approximate linear forward model, two algorithms based on the concept of Bayesian A-optimality were introduced. The first one can be run offline to find informative positions for the available electrodes assuming a Gaussian prior for the conductivity perturbation and assigning a ROI inside the patient's head. The second algorithm is an extension of the first one, with the Gaussian prior for the optimization process produced by a reconstruction method based on an initial electrode configuration and a combination of sequential linearizations and the lagged diffusivity iteration. Both algorithms produced intuitively acceptable results in the sense that the electrodes moved closer to the area of interest during the optimization process. In particular, the first algorithm led to significant reductions in the expected squared $L^2$ reconstruction errors under the linearized forward model.

Interesting topics for future research include simultaneous optimization of the electrode positions and employed current patterns (cf.~\cite{Kaipio04,Kaipio07}) as well as studying the effect of the particular anatomy of the patient on the optimal electrode positions. In this work, we mainly ignored the choice of the feeding electrode, although it certainly affects optimal configurations and the associated expected reconstruction errors, and we only considered a single head anatomy. Moreover, the introduced gradient-based optimization scheme should in the future be compared with a sparsification approach \cite{alexanderian2014optimal,alexanderian2016fast,Haber10} for choosing the optimal electrode positions.

\section*{Acknowledgements}
We would like to thank Valentina Candiani for her help with the computational head model.

\bibliographystyle{acm}
\bibliography{headelopt-refs15.bib}

\begin{thebibliography}{10}

\bibitem{alexanderian2021optimal_review}
{\sc Alexanderian, A.}
\newblock Optimal experimental design for infinite-dimensional {B}ayesian
  inverse problems governed by {PDE}s: A review.
\newblock {\em Inverse Problems\/} (2021), 043001.

\bibitem{alexanderian2016bayesian}
{\sc Alexanderian, A., Gloor, P.~J., Ghattas, O., et~al.}
\newblock On {B}ayesian {A}- and {D}-optimal experimental designs in infinite
  dimensions.
\newblock {\em Bayesian Anal. 11}, 3 (2016), 671--695.

\bibitem{alexanderian2014optimal}
{\sc Alexanderian, A., Petra, N., Stadler, G., and Ghattas, O.}
\newblock A-optimal design of experiments for infinite-dimensional {B}ayesian
  linear inverse problems with regularized {$l_0$}-sparsification.
\newblock {\em SIAM J. Sci. Comput. 36}, 5 (2014), A2122--A2148.

\bibitem{alexanderian2016fast}
{\sc Alexanderian, A., Petra, N., Stadler, G., and Ghattas, O.}
\newblock A fast and scalable method for {A}-optimal design of experiments for
  infinite-dimensional {B}ayesian nonlinear inverse problems.
\newblock {\em SIAM J. Sci. Comput. 38}, 1 (2016), A243--A272.

\bibitem{Bardsley18}
{\sc Bardsley, J.~M.}
\newblock {\em Computational uncertainty quantification for inverse problems},
  vol.~19 of {\em Computational Science \& Engineering}.
\newblock Society for Industrial and Applied Mathematics (SIAM), Philadelphia,
  PA, 2018.

\bibitem{Borcea02}
{\sc Borcea, L.}
\newblock Electrical impedance tomography.
\newblock {\em Inverse problems 18\/} (2002), R99--R136.

\bibitem{Burger21}
{\sc Burger, M., Hauptmann, A., Helin, T., Hyv\"onen, N., and Puska, J.-P.}
\newblock Sequentially optimized projections in x-ray imaging.
\newblock {\em Inverse Problems 37\/} (2014), 0750006.

\bibitem{Calvetti08}
{\sc Calvetti, D., and Somersalo, E.}
\newblock Hypermodels in the {B}ayesian imaging framework.
\newblock {\em Inverse Problems 24\/} (2008), 034013.

\bibitem{Candiani19}
{\sc Candiani, V., Hannukainen, A., and Hyv{\"o}nen, N.}
\newblock Computational framework for applying electrical impedance tomography
  to head imaging.
\newblock {\em SIAM J. Sci. Comput. 41\/} (2019), B1034--B1060.

\bibitem{Candiani21}
{\sc Candiani, V., Hyvönen, N., Kaipio, J.~P., and Kolehmainen, V.}
\newblock Approximation error method for imaging the human head by electrical
  impedance tomography.
\newblock {\em Inverse Problems 37}, 12 (2021), 125008.

\bibitem{Candiani20}
{\sc Candiani, V., and Santacesaria, M.}
\newblock Neural networks for classification of stroke in electrical impedance
  tomography on a 3{D} head model.
\newblock {\em Mathematics in Engineering 4\/} (2022), 1--22.

\bibitem{chaloner1995bayesian}
{\sc Chaloner, K., and Verdinelli, I.}
\newblock Bayesian experimental design: {A} review.
\newblock {\em Stat. Sci.\/} (1995), 273--304.

\bibitem{chan1999convergence}
{\sc Chan, T.~F., and Mulet, P.}
\newblock On the convergence of the lagged diffusivity fixed point method in
  total variation image restoration.
\newblock {\em SIAM J. Numer. Anal. 36}, 2 (1999), 354--367.

\bibitem{Cheney99}
{\sc Cheney, M., Isaacson, D., and Newell, J.}
\newblock Electrical impedance tomography.
\newblock {\em SIAM Rev. 41\/} (1999), 85--101.

\bibitem{Cheng89}
{\sc Cheng, K.-S., Isaacson, D., Newell, J.~S., and Gisser, D.~G.}
\newblock Electrode models for electric current computed tomography.
\newblock {\em IEEE Trans. Biomed. Eng. 36\/} (1989), 918--924.

\bibitem{Colton98}
{\sc Colton, D., and Kress, R.}
\newblock {\em Inverse acoustic and electromagnetic scattering theory},
  second~ed., vol.~93 of {\em Applied Mathematical Sciences}.
\newblock Springer-Verlag, Berlin, 1998.

\bibitem{Darde12}
{\sc Dard\'e, J., Hakula, H., Hyv\"onen, N., and Staboulis, S.}
\newblock Fine-tuning electrode information in electrical impedance tomography.
\newblock {\em Inverse Probl. Imag. 6\/} (2012), 399--421.

\bibitem{Darde21}
{\sc Dard{\'e}, J., Hyvönen, N., Kuutela, T., and Valkonen, T.}
\newblock Contact adapting electrode model for electrical impedance tomography.
\newblock {\em SIAM J. Appl. Math. 82\/} (2022), 427--229.

\bibitem{dobson1997convergence}
{\sc Dobson, D.~C., and Vogel, C.~R.}
\newblock Convergence of an iterative method for total variation denoising.
\newblock {\em SIAM J. Numer. Anal. 34}, 5 (1997), 1779--1791.

\bibitem{duong2022stability}
{\sc Duong, D.-L., Helin, T., and Rojo-Garcia, J.~R.}
\newblock Stability estimates for the expected utility in {B}ayesian optimal
  experimental design.
\newblock {\em Inverse Problems 39\/} (2022), 125008.

\bibitem{Gabriel96}
{\sc Gabriel, S., Lau, R.~W., and Gabriel, C.}
\newblock The dielectric properties of biological tissues: {II}. {M}easurements
  in the frequency range 10 {H}z to 20 {GH}z.
\newblock {\em Phys. Med. Biol. 41\/} (1996), 2251.

\bibitem{Hanke11b}
{\sc Hanke, M., Harrach, B., and Hyv{\"o}nen, N.}
\newblock Justification of point electrode models in electrical impedance
  tomography.
\newblock {\em Math. Models Methods Appl. Sci. 21\/} (2011), 1395--1413.

\bibitem{Hannukainen20}
{\sc Hannukainen, A., Hyv{\"o}nen, N., and Perkki{\"o}, L.}
\newblock Inverse heat source problem and experimental design for determining
  iron loss distribution.
\newblock {\em SIAM J. Sci. Comput 43\/} (2021), B243--B270.

\bibitem{Harhanen15}
{\sc Harhanen, L., Hyv\"onen, N., Majander, H., and Staboulis, S.}
\newblock Edge-enhancing reconstruction algorithm for three-dimensional
  electrical impedance tomography.
\newblock {\em SIAM J. Sci. Comput. 37\/} (2015), B60--B78.

\bibitem{Helin23}
{\sc Helin, T., Hyv\"onen, N., J., M., and Puska, J.-P.}
\newblock Bayesian design of measurements for magnetorelaxometry imaging.
\newblock {\em Inverse Problems\/} (2023), 125020.

\bibitem{Helin22}
{\sc Helin, T., Hyv\"onen, N., and Puska, J.-P.}
\newblock Edge-promoting adaptive {B}ayesian experimental design for x-ray
  imaging.
\newblock {\em SIAM J. Sci. Comput. 44\/} (2022), B506--B530.

\bibitem{Haber10}
{\sc Horesh, L., Haber, E., and Tenorio, L.}
\newblock {\em {\em Optimal Experimental Design for the Large-Scale Nonlinear
  Ill-Posed Problem of Impedance Imaging}, in `Large-Scale Inverse Problems and
  Quantification of Uncertainty'}.
\newblock Wiley, 2010.

\bibitem{Hyvonen17b}
{\sc Hyv\"onen, N., and Mustonen, L.}
\newblock Smoothened complete electrode model.
\newblock {\em SIAM J. Appl. Math. 77\/} (2017), 2250--2271.

\bibitem{Hyvonen18}
{\sc Hyv\"onen, N., and Mustonen, L.}
\newblock Generalized linearization techniques in electrical impedance
  tomography.
\newblock {\em Numer. Math. 140\/} (2018), 95--120.

\bibitem{Hyvonen14}
{\sc Hyv\"onen, N., Sepp\"anen, A., and Staboulis, S.}
\newblock Optimizing electrode positions in electrical impedance tomography.
\newblock {\em SIAM J. Appl. Math. 74\/} (2014), 1831--1851.

\bibitem{Kaipio06}
{\sc Kaipio, J., and Somersalo, E.}
\newblock {\em Statistical and Computational Inverse Problems}, vol.~160.
\newblock Springer Science \& Business Media, 2006.

\bibitem{Kaipio04}
{\sc Kaipio, J.~P., Sepp\"anen, A., Somersalo, E., and Haario, H.}
\newblock Posterior covariance related optimal current patterns in electrical
  impedance tomography.
\newblock {\em Inverse Problems 20\/} (2004), 919--936.

\bibitem{Kaipio07}
{\sc Kaipio, J.~P., Sepp\"anen, A., Voutilainen, A., and Haario, H.}
\newblock Optimal current patterns in dynamical electrical impedance tomography
  imaging.
\newblock {\em Inverse Problems 23\/} (2007), 1201--1214.

\bibitem{Karimi20}
{\sc Karimi, A., Taghizadeh, L., and Heitzinger, C.}
\newblock Optimal {B}ayesian experimental design for electrical impedance
  tomography in medical imaging.
\newblock {\em Comput. Methods Appl. Mech. Eng. 373\/} (2021), 113489.

\bibitem{Lai05}
{\sc Lai, Y., Van~Drongelen, W., Ding, L., Hecox, K.~E., Towle, V.~L., Frim,
  D.~M., and He, B.}
\newblock Estimation of in vivo human brain-to-skull conductivity ratio from
  simultaneous extra-and intra-cranial electrical potential recordings.
\newblock {\em Clin. Neurophysiol. 116\/} (2005), 456--465.

\bibitem{Lassas04}
{\sc Lassas, M., and Siltanen, S.}
\newblock Can one use total variation prior for edge-preserving bayesian
  inversion?
\newblock {\em Inverse Problems 20}, 5 (2004), 1537.

\bibitem{Latikka01}
{\sc Latikka, J.~A., Hyttinen, J.~A., Kuurne, T.~A., Eskola, H.~J., and
  Malmivuo, J.~A.}
\newblock The conductivity of brain tissue: {C}omparison of results in vivo and
  in vitro measurement.
\newblock In {\em Proceedings of the 23rd Annual International Conference of
  the IEEE Engineering in Medicine and Biology Society\/} (Instanbul, Turkey,
  October 2001), vol.~4, IEEE, pp.~910--912.

\bibitem{Lee16}
{\sc Lee, E., Duffy, W., Hadimani, R., Waris, M., Siddiqui, W., Islam, F.,
  Rajamani, M., Nathan, R., and Jiles, D.}
\newblock Investigational effect of brain-scalp distance on the efficacy of
  transcranial magnetic stimulation treatment in depression.
\newblock {\em IEEE Trans. Magn. 52\/} (2016), 1--4.

\bibitem{McCann19}
{\sc McCann, H., Pisano, G., and Beltrachini, L.}
\newblock Variation in reported human head tissue electrical conductivity
  values.
\newblock {\em Brain Topogr. 32\/} (2019), 825--858.

\bibitem{Nocedal06}
{\sc Nocedal, J., and Wright, S.~J.}
\newblock {\em Numerical optimization}, second~ed.
\newblock Springer, New York, 2006.

\bibitem{Oostendorp00}
{\sc Oostendorp, T.~F., Delbeke, J., and Stegeman, D.~F.}
\newblock The conductivity of the human skull: results of in vivo and in vitro
  measurements.
\newblock {\em IEEE Trans. Biomed. Eng. 47\/} (2000), 1487--1492.

\bibitem{rainforth2023modern}
{\sc Rainforth, T., Foster, A., Ivanova, D.~R., and Smith, F.~B.}
\newblock Modern {B}ayesian experimental design.
\newblock {\em arXiv preprint arXiv:2302.14545\/} (2023).

\bibitem{Rudin92}
{\sc Rudin, L.~I., Osher, S., and Fatemi, E.}
\newblock Nonlinear total variation based noise removal algorithms.
\newblock {\em Physica D 60\/} (1992), 259--268.

\bibitem{ryan2016review}
{\sc Ryan, E.~G., Drovandi, C.~C., McGree, J.~M., and Pettitt, A.~N.}
\newblock A review of modern computational algorithms for {B}ayesian optimal
  design.
\newblock {\em Int. Stat. Rev. 84}, 1 (2016), 128--154.

\bibitem{Smyl20}
{\sc Smyl, D., and Liu, D.}
\newblock Optimizing electrode positions in 2-{D} electrical impedance
  tomography using deep learning.
\newblock {\em {IEEE} Trans. Instrum. Meas. 69}, 9 (2020), 6030--6044.

\bibitem{Somersalo92}
{\sc Somersalo, E., Cheney, M., and Isaacson, D.}
\newblock Existence and uniqueness for electrode models for electric current
  computed tomography.
\newblock {\em SIAM J. Appl. Math. 52\/} (1992), 1023--1040.

\bibitem{Toivanen21}
{\sc Toivanen, J., H\"anninen, A., Savolainen, T., Forss, N., and Kolehmainen,
  V.}
\newblock 8 - monitoring hemorrhagic strokes using eit.
\newblock In {\em Bioimpedance and Spectroscopy}, P.~Annus and M.~Min, Eds.
  Academic Press, 2021, pp.~271--298.

\bibitem{Uhlmann09}
{\sc Uhlmann, G.}
\newblock Electrical impedance tomography and {C}alder{\'o}n's problem.
\newblock {\em Inverse Problems 25\/} (2009), 123011.

\bibitem{Vauhkonen97}
{\sc Vauhkonen, M.}
\newblock {\em Electrical impedance tomography with prior information},
  vol.~62.
\newblock Kuopio University Publications C (Dissertation), 1997.

\bibitem{Vauhkonen99}
{\sc Vauhkonen, P.~J., Vauhkonen, M., Savolainen, T., and Kaipio, J.~P.}
\newblock Three-dimensional electrical impedance tomography based on the
  complete electrode model.
\newblock {\em IEEE Trans. Biomed. Eng. 46\/} (1999), 1150–1160.

\bibitem{Vogel96}
{\sc Vogel, C.~R., and Oman, M.~E.}
\newblock Iterative methods for total variation denoising.
\newblock {\em SIAM J. Sci. Comput. 17\/} (1996), 227–238.

\end{thebibliography}
\end{document}